\magnification = 1200
\parskip = 6pt
\input pictex.tex

\font\title=cmbx12

\def \ZZ{\hbox {Z\kern -3pt Z}}
\def \RR{\hbox{I\kern -2pt R}}
\def \PP{\hbox{I\kern -2pt P}}
\def \TT{\hbox{I\kern -4pt T}}
\def\FF{\hbox {I\kern -2pt F}}
\def\NN{\hbox {I\kern -2pt N}}
\def\QQ{\hbox {I\kern -2pt Q}}
\def\CC{\hbox {I\kern -2pt C}}
\def\SS{\hbox {I\kern -2pt S}}
\def\FF{\hbox {I\kern -2pt F}}
\def\NN{\hbox {I\kern -2pt N}}
\def\QQ{\hbox {I\kern -2pt Q}}
\def\CC{\hbox {I\kern -2pt C}}
\def\SS{\hbox {I\kern -2pt S}}
\def\hollowbox{\hbox{\vbox{\hrule\hbox{\vrule\vbox{\vbox to 8pt{\vfill\hbox to
5pt{\hfill}}}\vrule}\hrule}}}
\def\ep{\hfill\break\null\hskip 15.5truecm \hollowbox\goodbreak\medskip}
\def \cross{C(\Omega) \rtimes \Gamma}
\def \crossr{C(\Omega) \rtimes_r \Gamma}
\def \G{\Gamma}
\def \O{\Omega}
\def \o{\omega}

\hfuzz = 50pt
\font\tenmsb=msbm10
\font\sevenmsb=msbm7
\font\fivemsb=msbm5
\newfam\msbfam
\textfont\msbfam=\tenmsb
\scriptfont\msbfam=\sevenmsb
\scriptscriptfont\msbfam=\fivemsb
\def\hexnumber@#1{\ifnum#1<10 \number#1\else
 \ifnum#1=10 A\else\ifnum#1=11 B\else\ifnum#1=12 C\else
 \ifnum#1=13 D\else\ifnum#1=14 E\else\ifnum#1=15 F\fi\fi\fi\fi\fi\fi\fi}
\def\msb@{\hexnumber@\msbfam}
\mathchardef\rtimes="2\msb@6F

\centerline {\bf   $C^*$--algebras arising from group actions  on the boundary of a triangle building}
\bigskip
\centerline{\it  Guyan Robertson\footnote *{This research was supported by the Australian
 Research Council}   and  Tim Steger*}

\bigskip
\medskip

\medskip
\baselineskip = 20pt.

\noindent 0. {\bf Introduction}.

A subgroup of an amenable group is amenable.  The $C^*$-algebra version of this fact is false.
This was first proved by M.-D. Choi [6] who proved that the non-nuclear
$C^*$-algebra $C^*_r(\ZZ_2*\ZZ_3)$ is a subalgebra of the nuclear  Cuntz algebra ${\cal O}_2$.
A. Connes provided another example , based on a crossed product construction. More recently
J. Spielberg [23] showed that these examples were essentially the same. In fact he proved
that certain of the $C^*$-algebras studied by J. Cuntz and W. Krieger [10] can be constructed
naturally as crossed product algebras. For example if the group $\Gamma$ acts simply transitively
on  a homogeneous tree of finite degree with boundary $\Omega$ then $\cross$
 is a Cuntz-Krieger algebra.

 Such trees may be regarded as affine buildings of type $\widetilde A_1$. The present paper is
devoted to the study of the analogous situation where a group $\G$ acts simply transitively
on the vertices of an affine building of type $\widetilde A_2$ with boundary $\O$ [7].
The corresponding crossed product algebra $\cross$ is then generated by two Cuntz-Krieger
 algebras . (See \S3.) Moreover we show that $\cross$ is simple and nuclear. This is a consequence
of the facts that the action of $\G$ on $\O$  is minimal, topologically free, and amenable.
( See \S4.)

\vfill\eject
\noindent 1. {\bf Background}.

The original motivation was an attempt to show that the reduced $C^*$--algebra of a
discrete group is subnuclear (i.e. a subalgebra of a nuclear $C^*$--algebra). M. D. Choi in
[6] gave the first result in this direction by embedding {$C^*_r(\ZZ_2 *\ZZ_3)$} in the
Cuntz algebra ${\cal O}_2$.  At about the same time A. Connes gave other examples in a talk at
Kingston in 1980.  The idea of Connes was included in the first version of [14] and was also
rediscovered by J. Quigg and J. Spielberg as part (a) of the following result [21:  Corollary 4.3].
See also [2: Th\'eor\`eme 4.6], and the remark following it.

\noindent {\bf Proposition 1.1}.  Let $G$ be a locally compact group and let $P$ and
$\Gamma$ be closed subgroups of $G$ such that $P$ is amenable.   Then

\noindent (a)  $C_0(G/P) \rtimes \Gamma$ is a nuclear $C^*$--algebra;

\noindent (b)   $C_0(G/P) \rtimes \Gamma$ is canonically isomorphic to the reduced
crossed product  $C_0(G/P) \rtimes_r \Gamma$.
\bigskip
\noindent {\bf Connes' Example}.  Take $\Gamma = PSL(2, \ZZ) \cong \ZZ_2
*\ZZ_3$.  Let $G = PSL(2, \RR)$ and $P$ the subgroup of upper triangular matrices.
Then $G$ acts transitively on $\PP_1\RR = \TT$ and $P = \hbox { Stab } (\infty)$.   Hence
$G/P \cong \TT$.  By Proposition 1, we have that $C(\TT) \rtimes_r \Gamma$ is nuclear.
Therefore we have an embedding of $C^*_r(\ZZ_2 *\ZZ_3)$ into a nuclear $C^*$--algebra
which at first sight appears to be different from Choi's embedding in ${\cal O}_2$.

Note that $G/P$ is the Furstenberg boundary of $G$ [12].  It turns out that the approaches
of Choi and Connes may be unified by using the concept of a boundary.
\bigskip
\bigskip

\noindent 2. {\bf Cuntz--Krieger algebras arising from boundary actions of free products of
cyclic groups}.

The ideas of this section are due to J. Spielberg [23].  They provide the motivation for our results
for groups acting on buildings.

We present the case $\Gamma = \ZZ * \ZZ$, the free group with generators $a$ and $b$.
The Cayley graph of $\Gamma$ is a homogeneous tree of degree 4.  The vertices of the tree
are elements of $\Gamma$ i.e. reduced words in the generators and their inverses.

The {\it boundary} $\Omega$ of the tree can be thought of as the set of all infinite reduced
words $\omega = x_1 x_2 x_3 \ldots$, where $x_i \in S = \left \{a, b, a^{-1}, b^{-1} \right \}$.
$\Omega$ has a natural compact (totally disconnected) topology in which an open
neighbourhood of $\omega \in \Omega$ consists of those $\omega^\prime \in
\Omega$ whose corresponding infinite word agrees with that of $\omega$ on a finite initial
segment [11].  Left multiplication by $x \in \Gamma$ defines a homeomorphism of
$\Omega$ and so induces an action $\alpha$ of $\Gamma$ by
$$
\alpha (x) f(\omega) = f(x^{-1} \omega).
$$

For $x \in \G$ let $\Omega (x)$ be the set of
infinite words beginning with $x$.  Then $\Omega(x)$ is open and closed in $\Omega$ and
the sets $g \Omega(x)$ and $g\left (\Omega \backslash
 \Omega(x)\right)$ , where $g \in \G$ and $x \in S$, form a base for
the topology of $\Omega$.

We partition the boundary $\Omega$ into four parts as follows.

\midinsert
  \line \bgroup \hss
\beginpicture
\setcoordinatesystem units <0.5cm,0.5cm>
\put {$ Figure \  2.1$} [t] at 0 -7.5
\put{$a$}   [b,r] at 1.8 0.2
\put{$b$}   [t,l] at 0.2 1.8
\put{$a^{-1}$}   [b,l]  at  -1.8  0.2
\put{$b^{-1}$}   [b,l]  at   0.2 -1.6
\put{$\Omega(a)$}  [l]  at  5.5 0
\put{$\Omega(b)$}  [b] at  0  5.5
\put{$\Omega(a^{-1})$}  [r] at  -5.5  0
\put{$\Omega(b^{-1})$}  [t] at  0  -5.5
\putrule from -4 0 to 4 0
\putrule from 0 -4 to 0 4
\putrule from 2 -1.5 to 2 1.5
\putrule from -1.5 2 to 1.5 2
\putrule from -2 -1.5 to -2 1.5
\putrule from -1.5 -2 to 1.5 -2
\putrule from 3.5 -.5 to 3.5 .5
\putrule from 1.5 1 to 2.5 1
\putrule from 1.5 -1 to 2.5 -1
\putrule from -.5 3.5 to .5 3.5
\putrule from  1 1.5 to 1 2.5
\putrule from  -1 1.5 to -1 2.5
\putrule from  -3.5 .5 to -3.5 -.5
\putrule from  -1.5 1 to -2.5 1
\putrule from  -1.5 -1 to -2.5 -1
\putrule from -.5 -3.5 to .5 -3.5
\putrule from  1  -1.5 to 1 -2.5
\putrule from  -1 -1.5 to -1 -2.5
\setdashes
\circulararc 80 degrees from 3.83 -3.214 center at 0 0
\circulararc 80 degrees from 3.214 3.83  center at 0 0
\circulararc 80 degrees from -3.83 3.214 center at 0 0
\circulararc 80 degrees from -3.214 -3.83 center at 0 0
   \endpicture
  \hss \egroup
\endinsert

Recall that the crossed product $C(\Omega) \rtimes \Gamma$ is the universal covariant
representation of $\left ( C(\Omega), \Gamma, \alpha \right )$.  In other words,
$C(\Omega) \rtimes \Gamma$ is generated by $C(\Omega)$ and the image of a unitary
representation $\pi$ of $\Gamma$ such that $\alpha (g) f = \pi(g) f \pi(g)^*$ for
$f \in C(\Omega)$ and $g \in \Gamma$, and every such $C^*$--algebra is a quotient
of $C(\Omega) \rtimes \Gamma$. It is convenient to simply write  $g$ instead of $\pi(g)$, thereby
identifying elements of $\Gamma$ with unitaries in $C(\Omega) \rtimes \Gamma$. Also
for $x \in \Gamma$ let $p_x$ denote the projection defined by the characteristic function
$1_{\Omega(x)} \in C(\Omega)$.

For $g \in \Gamma$, we have
$$
g\  p_x g^{-1} = \alpha(g) 1_{\Omega(x)} = 1_{g \Omega (x)}
$$
and therefore, for each $x \in S$,
$$
p_x + x\  p_{x^{-1}}\ x^{-1} = 1.
$$
Moreover
$$
p_a + p_{a^{-1}} + p_b + p_{b^{-1}} = 1
$$
and covariant representations of $\left (C(\Omega), \Gamma, \alpha\right )$ correspond
to representations of these relations.

For $x \in S$, define a partial isometry $s_x \in C(\Omega) \rtimes \Gamma$ by
$$\eqalignno{
s_x & = x (1 - p_{x^{-1}}) \cr
\noalign {\hbox {Then}}
s_x s_x^* & = x(1 - p_{x^{-1}}) x^{-1} = p_x \cr
\noalign {\hbox { and }}
s_x^* s_x & = 1 - p_{x^{-1}} = \sum_{y \not= x^{-1}} s_y s_y^*  \cr}
$$
These relations show that the partial isometries $s_x, x \in S$, generate the Cuntz--Krieger
algebra ${\cal O}_A$ of [10], where
$$
A = \left [ \matrix {1&0&1&1 \cr
0&1&1&1 \cr
1&1&1&0 \cr
1&1&0&1 \cr} \right]
$$
relative to $\{a, a^{-1}, b, b^{-1}\} \times \{a, a^{-1}, b, b^{-1} \}$.

We can recover the generators for $C(\Omega) \rtimes \Gamma$ from the
$s_x$'s by
$$ \eqalignno {
x & = s_x + s_{x^{-1}}^* \cr
\noalign {\hbox {and }}
p_x & = s_x s_x^* .\cr}
$$
\noindent Therefore $C(\Omega) \rtimes \Gamma = {\cal O}_A$, which is a simple
$C^*$--algebra since $A$ is an irreducible matrix [10].  It follows that the reduced crossed
product  $C(\Omega) \rtimes_r \Gamma$, being a quotient of  $C(\Omega) \rtimes
\Gamma$  is also isomorphic to ${\cal O}_A$.  The algebra ${\cal O}_A$ is nuclear [10]
 and so we obtain in
particular that $C_r^* (\Gamma)$ is subnuclear.

In [23], J. Spielberg covers the case where $\Gamma$ is a free product of cyclic groups of
arbitrary order.  In particular this method gives an alternative demonstration of Choi's
embedding of $C_r^* (\ZZ_2 * \ZZ_3)$ into ${\cal O}_2$.

\noindent {\bf Remarks.}(1) The fact that the algebras ${\cal O}_A$ obtained in this way are
nuclear may also be seen by using the fact that the action of $\Gamma$ on $\Omega$ is amenable in the
sense of [2]. This can be done by constructing a sequence of functions $(f_n)$ on
$\Gamma \times \Omega$ satisfying condition $(h)$ of [2: Th\'eor\`eme 4.9].. For the case
where $\Gamma$ is the free group on two generators this is done explicitly in the appendix
of [16]. See also [15].  More generally the arguments of [1] show how to do it for any hyperbolic group
$\Gamma$. It follows from [2: Th\'eor\`eme 4.5 and Proposition 4.8]  that the crossed
product algebras  $C(\Omega) \rtimes \Gamma$ and $C(\Omega) \rtimes_r \Gamma$  are canonically
isomorphic and nuclear.

(2) The definition of the operator $s_x$ may be expressed as follows.
$$s_x = x \sum_{y \in A_x}p_y$$
where  $A_x$ is the set of elements $y$ of the generating set  $S$  such that  $xy \ne e$. That is
 the Cayley graph of $\Gamma$ contains a semi-infinite geodesic of the form
\bigskip

\centerline{
\beginpicture
\setcoordinatesystem units <1cm,1cm>
\put {$ Figure\  2.2$} [t] at 1 -1
\multiput {$\bullet$} at -2 0 *3 2 0 /
\putrule from -2 0 to 2 0
\setdashes
\putrule from 2 0 to 5 0
\put {$e$}  [t]  at  -2 -0.2
\put {$x$}  [t]   at 0  -0.2
\put {$xy$}  [t]  at  2 -0.2
\endpicture
}
\bigskip
In this  geometric form, we shall extend the definitions to the case of affine buildings.
\bigskip

\noindent 3. {\bf Crossed product algebras associated with triangle buildings.}

The construction of the previous section shows that if one has a group $\Gamma$ acting
simply transitively on the vertices of a tree with boundary $\Omega$ then, relative to the
natural action on $\Omega$, the crossed product algebra $C(\Omega) \rtimes \Gamma$
can be identified as a (simple, nuclear) Cuntz--Krieger algebra.

We now give an analogous construction for groups acting on the vertices of affine buildings
of type $\widetilde A_2$.

\noindent 3.1 {\bf \sl Triangle Buildings}

A triangle building is a thick affine building $\triangle$ of type  $\widetilde A_2$ [22].
This means that $\triangle$ is a chamber system consisting of vertices, edges and triangles
({\sl chambers}).  Each edge lies on $q + 1$ triangles, where $q \geq 2$ is the {\sl order} of
$\triangle$.  An {\sl apartment} is a subcomplex of $\triangle$ isomorphic to the Euclidean
plane tesselated by equilateral triangles.  A {\sl sector} (or {\sl Weyl chamber}) is a
 sector in some apartment

\midinsert
 \line \bgroup \hss
\beginpicture
\setcoordinatesystem units <0.5cm,0.866cm>   
\setplotarea x from -4 to 4, y from  -1 to  4  
\put {$ Figure \  3.1$} [t] at 0  -1
\putrule from -1 1 to 1 1
\putrule from -2 2 to 2 2
\setlinear \plot -2 2  0 0  2 2  /
\setlinear \plot -1 1  0 2  1 1 /
\setdashes
\setlinear \plot -3.5 3.5  -2 2  -1 3  0 2  1 3  2 2  3.5 3.5 /
\putrule from  -3 3  to  3 3
\endpicture
 \hss \egroup
\endinsert

\noindent Two sectors are  {\sl equivalent} (or parallel) if their intersection contains a sector.
In the case of a tree, sectors correspond to semi-infinite geodesics and apartments to (doubly
infinite) geodesics.

\noindent 3.2 {\bf \sl The boundary $\Omega$ of $\triangle$}

The boundary $\Omega$ is defined to be the set of equivalence classes of sectors in $\triangle$.  Fix a
vertex $x$.  For any $\omega \in \Omega$ there is a unique sector $S^x(\omega)$ in the
class $\omega$ having base vertex $x$ [22: Theorem 9.6].
\bigskip

\midinsert                        
 \line \bgroup \hss
\beginpicture
\setcoordinatesystem units <0.5cm,0.866cm>   
\put {$Figure \ 3.2.1$} at 0 -1.5
\put {$\bullet$} at 0 0
\multiput {$\bullet$} at -1 1 *1 2 0 /
\multiput {$\bullet$} at -2 2 *2 2 0 /
\put{$x$}     [t]   at  0  -0.2
\put{$\omega_{1,1}$} [b r] at  -0.2 2.05
\put{$\omega_{0,1}$} [r]  at -1.2 1
\put{$\omega_{0,2}$}  [r]  at -2.2 2
\put{$\omega_{1,0}$}  [l] at  1.2 1
\put{$\omega_{2,0}$}  [l] at 2.2 2
\putrule from -1 1 to 1 1
\putrule from -2 2 to 2 2
\setlinear \plot -2 2  0 0  2 2  /
\setlinear \plot -1 1  0 2  1 1 /
\setdashes
\setlinear \plot -3.5 3.5  -2 2  -1 3  0 2  1 3  2 2  3.5 3.5 /
\putrule from  -3 3  to  3 3
\endpicture
  \hss \egroup
\endinsert

\bigskip
In the terminology of [5: Chapter VI.9]
$\Omega$ is the set of chambers of the building at infinity $\triangle^\infty$.
Also $\Omega$ is a totally disconnected compact Hausdorff space with a base for the
topology being given by sets of the form
$$
\Omega^x(v) = \left \{ \omega \in \Omega : S^x(\omega) \hbox { contains } v \right \}
$$
where $v$ is a vertex of $\triangle$ [7 : \S 2].
\bigskip

\midinsert                                        
  \line \bgroup \hss
\beginpicture
\setcoordinatesystem units <0.5cm, 0.866cm>
\put {$ Figure \ 3.2.2$} at 0 -1
\put {$\bullet$} at 0 0
\put {$\bullet$} at -1 2
\put {$S^{x}(\omega)$}  at 0 4
\put {$x$} [t] at 0 -0.2
\put {$v$}  at  -0.8 2.2
\setlinear \plot -4 4  0 0  4 4 /
\endpicture
  \hss \egroup
\endinsert
\bigskip

\vfill \eject

\noindent 3.3 {\bf \sl The group $\Gamma$}

The buildings we consider are constructed in [8].  They are exactly the triangle
buildings on whose vertices a group acts simply transitively and in a ``type rotating'' way.

Let $(P,L)$ be a projective plane of order $q$.  There should be no confusion between this
$P$ and the group $P$ of section 1.  There are $q^2 + q + 1$ points (elements of $P$) and $q^2+q+1$
lines (elements of $L$).  Each point lies on $q+1$ lines and each line contains $q+1$ points.
Let $\lambda : P \rightarrow L$ be a bijection (a point--line correspondence).  Let ${\cal T}$
be a set of triples $(x, y, z)$ where $x, y, z \in P$, with the following properties.

\item {(i)}  Given $x, y \in P$, then $(x, y, z) \in {\cal T}$ for some
$z \in P$ if and only if $y$ and $\lambda(x)$ are incident (i.e. $y \in
\lambda(x)$).

\item {(ii)}  $(x, y, z) \in {\cal T} \Rightarrow (y, z, x) \in {\cal T}$.

\item {(iii)}  Given $x, y \in P$, then $(x, y, z) \in {\cal T}$ for at
most one $z \in P$.

${\cal T}$ is called a {\it triangle presentation} (or {\it triella}) compatible with
$\lambda$.  A complete list is given in [8] of all triangle presentations for $q = 2$
and $q = 3$.

Let $\{a_x : x \in P\}$ be  $q^2 + q + 1$ distinct letters and form the group
$$
\Gamma = \big\langle a_x, x \in P \ |\  a_x a_y a_z = 1 \hbox { for } (x, y, z) \in {\cal T}
\big \rangle $$
\noindent Then ${\cal T}$ gives rise to a triangle building $\triangle_{\cal T}$ whose
vertices and edges form the Cayley graph of $\Gamma$ with respect to the
generators $a_x, x \in P$, and their inverses, and whose chambers are the sets
$$
\{ g, g a_x^{-1}, g a_y\}
$$
where $g \in \Gamma$ and $(x, y, z) \in {\cal T}$, for some $z$.
\bigskip

\centerline{
\beginpicture
\setcoordinatesystem units <1cm, 1.732cm>
\put {$Figure\ 3.3.1$} [t] at 0 -0.5
\put {$\bullet$} at 0 0
\multiput {$\bullet$} at -1 1 *1 2 0 /
\put {$g$} [t] at 0 -0.1
\put {$ga_x^{-1}$} [r] at -1.1 1
\put {$ga_y$} [l] at 1.1 1
\putrule from -1 1 to 1 1
\setlinear \plot -1 1 0 0 1 1 /
\endpicture
}
\bigskip

\noindent This chamber can also be represented by the diagram

\midinsert
 \line \bgroup \hss
\beginpicture
\setcoordinatesystem units  <1cm, 1.732cm> point at 5 0
\put {$Figure \ 3.3.2$} [t] at  0 -0.5
\put {$\bullet$} at 0 0
\multiput {$\bullet$} at -1 1 *1 2 0 /
\put {$g$} [t] at 0 -0.1
\arrow <10pt> [.2, .67] from  0.2 1 to 0 1
\arrow <10pt> [.2, .67] from  -0.7 0.7 to -0.5 0.5
\arrow <10pt> [.2, .67] from  0.3 0.3 to 0.5 0.5
\put {$a_x$} [r ] at -0.6 0.5
\put {$a_y$} [ l] at 0.6 0.5
\put {$a_z$} [ b] at 0 1.1
\putrule from -1 1 to 1 1
\setlinear \plot -1 1 0 0 1 1 /
\endpicture
 \hss \egroup
\endinsert

\bigskip
It is often convenient to identify the point  $x \in P$ with the generator $a_x \in \G$
The lines in $L$ correspond to the inverse generators of $\G$ and the
point-line correspondence $\lambda : P \rightarrow L$ satisfies
$a_{\lambda(x)} = a_x^{-1}$ for $x \in P$ [8]. We may therefore write
$x^{-1}$ for $a_x^{-1}$ and identify $x^{-1}$ with  $\lambda(x)$.

The set of nearest neighbours of an element $g \in \G$  (the {\it residue} of $g$)
can be identified with a projective plane  of order $q$, having  $q^2 + q + 1$ points
$\{ gx : x \in P \}$ and $q^2+q+1$ lines $\{ gx^{-1} : x \in P \}$.  According to
condition (1) above the line  $gx^{-1}$ is incident with the point  $gy$
if and only if  $xyz = 1$ for some $z \in P$. This is equivalent to
$gx^{-1} = gyz$, which means that $gx^{-1}$ and  $gy$ are adjacent in the Cayley
graph of $\G$, as in Figure 3.3.1.

From now on let $\triangle$ be a triangle building arising from a triangle presentation
as above and let
$\Gamma$ be the corresponding group which acts simply transitively on the vertices.  We
may identify $\Gamma$ with the vertices of its Cayley graph i.e. with the vertices of
$\triangle$.

 Any element
$g \in \Gamma \backslash \{ e \}$ can be written uniquely in the left normal form
$$g = x_1^{-1}x_2^{-1} \dots x_n^{-1}y_1y_2 \dots y_m$$
where there are no obvious cancellations and  $x_i,y_j \in P$, $1 \le i \le n, 1 \le j \le m$.
(As usual, the concise notation $x,y$ has been used in place of $a_x,a_y$.) [9: Lemma 6.2].
The absence of ``obvious'' cancellations means that
$ x_i \notin \lambda(x_{i+1})  \quad ( 1 \le i < n)$ ,
$ y_{j+1} \notin \lambda(y_j)  \quad (1 \le j < m)$, and
$ x_n \ne y_1 $.
Also any such word for $g$ is a minimal word for $g$ in the generators $x \in P$ and their inverses
[9:Lemma 6.2].

It follows from [9: Lemma 6.2] and the remarks following it, that if a group
 element $g$ is expressed as above in left normal form
then its position in a sector based at $e$ containing it is obtained by moving $n$ steps up
the ``left hand'' wall of the sector then $m$ steps in the direction of the ``right hand'' wall.
( The choice of ``left'' and ``right'' is purely conventional in our digrams.)  Moreover we can
fill out any apartment containing $e$ by starting at $e$ and repeatedly multiplying on the right by
a generator or inverse generator in a way consistent with the triangle presentation of the
group. An apartment is a union of six sectors based at $e$.

\midinsert
  \line \bgroup \hss
\beginpicture
\setcoordinatesystem units <1cm, 1.732cm>
\setplotarea  x from -2 to 2,  y from -1.5 to 1
\put {$\it Figure \ 3.3.3$} at  0 -1.5
\putrule from -2 0 to 2 0
\setlinear \plot  -1 -1  0 0  1 1  /
\setlinear \plot  1 -1  0 0  -1 1 /
\endpicture
  \hss \egroup
\endinsert

A simple example arises when $q = 2$ and $\Gamma$ is taken to have
 generators $a_0, a_1, \ldots, a_6$
and relations
$$
a_{[i]} \ a_{[i+1]} \ a_{[i+3]} = 1,  \quad \quad i\in \{ 0, \ldots , 6\}
$$
where $[j]$ denotes $j$ mod 7.  See [8 : Part II, $\S$4].
Then Figure 3.3.4 represents part of an apartment  which contains
$e$.

\midinsert
  \line \bgroup \hss
\beginpicture
\setcoordinatesystem units  <1cm, 1.732cm>        
\setplotarea x from -2.5 to 5, y from -2 to 2    
\put {$Figure \  3.3.4$} at  0 -2
\put {$a_3^{-2}$} [b] at -2  2.1
\put {$a_3^{-1}a_5$} [b] at  0  2.1
\put {$a_0a_3$}    [b] at  2  2.1
\put {$a_6^{-1}a_4$} [l] at 3.1 1
\put {$a_6^{-1}a_0^{-1}$} [l] at 4.1 0
\put {$a_6^{-1}a_1$} [t] at 3  -1.1
\put {$a_3$}  [t] at 1  -1.2
\put {$a_1^{-1}$} [t] at  -1  -1.1
\put {$a_1^{-1}a_3$} [t] at -3 -1.1
\put {$a_5^2$} [r]  at  -4.1 0
\put {$a_3^{-1}a_4$} [r] at -3.2  1
\put {$a_3^{-1}$} [l,b] at -0.7 1.1
\put {$a_0$}  [l,b] at 1.3  1.1
\put {$a_6^{-1}$} [l] at  2.3 0.2
\put {$a_5$}  [l,b]  at  -1.7 0.1
\put {$e$} [b]  at 0.3  0.1
\putrule from -2   2     to  2  2
\putrule from  -3 1  to 3 1
\putrule from -4 0  to 4 0
\putrule from -3 -1  to  3 -1
\setlinear
\plot  -4 0  -2 2 /
\plot -3 -1  0 2 /
\plot -1 -1  2 2 /
\plot  1 -1  3 1 /
\plot  3 -1  4 0 /
\plot  -4 0  -3 -1 /
\plot  -3 1   -1 -1 /
\plot  -2 2   1 -1 /
\plot  0 2   3 -1 /
\plot  2 2   4 0 /
\endpicture
  \hss \egroup
\endinsert

For the purposes of constructing such apartments it is more convenient to label
the edges as directed line segments, where the label corresponds to right multiplication
by a generator or an inverse generator. For example,the initial portion of
 one of the six sectors
based at $e$ in this apartment may be represented as in Figure 3.3.5.

\midinsert
  \line \bgroup \hss
\beginpicture
\setcoordinatesystem units <1cm, 1.732cm>
\setplotarea  x from -2 to 2,  y from -2 to 2
\put {$Figure \  3.3.5$} at  0 -1.5
\put {$a_3$} [r ] at -0.6 -0.5
\put {$a_0$} [ l] at 0.6 -0.5
\put {$a_6$} [r ] at 0.4  0.5
\put {$a_3$} [ l] at 1.6  0.5
\put {$a_3$} [r ] at -1.5  0.4
\put {$a_5$} [ l] at -0.4  0.5
\put {$a_1$} [t] at  0.0   -0.1
\put {$a_2$} [b] at  -1.0  0.8
\put {$a_4$} [b] at  1.0   0.8
\put {$e$} [t] at  0.0  -1.1
\multiput {\beginpicture
\setcoordinatesystem units <1cm, 1.732cm>
\arrow <10pt> [.2, .67] from  0.2 1 to 0 1
\arrow <10pt> [.2, .67] from  -0.7 0.7 to -0.5 0.5
\arrow <10pt> [.2, .67] from  0.3 0.3 to 0.5 0.5
\putrule from -1 1 to 1 1
\setlinear \plot -1 1 0 0 1 1 /
\endpicture}  at   0 0  1 1  -1 1  /
\endpicture
  \hss \egroup
\endinsert

Since $q = 2$ in this example, the projective plane of nearest neighbours of
$e$ has  seven  points  and seven lines. Thus the residue of $e$ consists
of fourteen vertices. A vertex representing a point is adjacent to a vertex representing
a line if the point and line are incident in the projective plane.
 The adjacencies between these vertices
are shown in Figure 3.3.6. (Every vertex shown is of course also adjacent to $e$.)

Each of the 28 hexagons contained in Figure 3.3.6 lies in a different
apartment of $\triangle$ containing $e$.

\bigskip
\bigskip

\midinsert
  \line \bgroup \hss
\beginpicture

\setcoordinatesystem units <3cm,3cm>
\setplotarea  x  from  -1.5  to  1.5, y  from  -1.5  to  1.5
\put {$\it Figure \  3.3.6$} at 0 -1.5

\put {$\bullet$} at 1 0

\put {$\bullet$} at 0.9009690323808488999   0.4338833975744138177

\put {$\bullet$} at 0.6234903946185663099   0.7818310097574695578

\put {$\bullet$} at 0.222522042695637753     0.9749276591186449509

\put {$\bullet$} at -0.2225194556367689247  0.9749282495974337219

\put {$\bullet$} at -0.6234883199575836131   0.7818326642424642541

\put {$\bullet$} at -0.9009678810291215377   0.433885788375114562

\put {$\bullet$} at -1 0

\put {$\bullet$} at 0.9009690323808488999   -0.4338833975744138177

\put {$\bullet$} at 0.6234903946185663099   -0.7818310097574695578

\put {$\bullet$} at 0.222522042695637753     -0.9749276591186449509

\put {$\bullet$} at -0.2225194556367689247  -0.9749282495974337219

\put {$\bullet$} at -0.6234883199575836131   -0.7818326642424642541

\put {$\bullet$} at -0.9009678810291215377   -0.433885788375114562

\put{$a_1$} [l] at 1.1 0

\put{$a_0^{-1}$} [l] at 1.0009690323808488999   0.4338833975744138177

\put{$a_2$} [l,b] at 0.6734903946185663099   0.8318310097574695578

\put{$a_1^{-1}$} [[b] at 0.222522042695637753     1.0749276591186449509

\put{$a_3$} [b] at -0.2225194556367689247  1.0749282495974337219

\put{$a_2^{-1}$} [b,r] at -0.6734883199575836131   0.8318326642424642541

\put{$a_4$} [r] at -1.0009678810291215377   0.433885788375114562

\put{$a_3^{-1}$} [r] at -1.1 0

\put{$a_6^{-1}$} [l] at 1.0009690323808488999   -0.4338833975744138177

\put{$a_0$} [l,t] at 0.6734903946185663099   -0.8318310097574695578

\put{$a_5^{-1}$} [t] at 0.222522042695637753     -1.0749276591186449509

\put{$a_6$} [t] at -0.2225194556367689247  -1.0749282495974337219

\put{$a_4^{-1}$} [t,r] at -0.6734883199575836131   -0.8318326642424642541

\put{$a_5$} [r] at -1.0009678810291215377   -0.433885788375114562

\setlinear \plot  1 0 
 0.9009690323808488999   0.433883397574413817    
 -0.9009678810291215377   0.433885788375114562   
 -1 0    
0.6234903946185663099   -0.7818310097574695578   
 0.9009690323808488999   -0.4338833975744138177  
 -0.2225194556367689247  0.9749282495974337219   
 -0.6234883199575836131   0.7818326642424642541  
-0.2225194556367689247  -0.9749282495974337219   
 0.222522042695637753     -0.9749276591186449509 
 0.6234903946185663099   0.7818310097574695578   
 0.222522042695637753     0.9749276591186449509  
 -0.9009678810291215377   -0.433885788375114562   
-0.6234883199575836131   -0.7818326642424642541  
1 0 
 0.9009690323808488999   -0.4338833975744138177  
0.6234903946185663099   -0.7818310097574695578   
0.222522042695637753     -0.9749276591186449509 
-0.2225194556367689247  -0.9749282495974337219   
-0.6234883199575836131   -0.7818326642424642541  
-0.9009678810291215377   -0.433885788375114562   
 -1 0    
 -0.9009678810291215377   0.433885788375114562   
-0.6234883199575836131   0.7818326642424642541  
-0.2225194556367689247  0.9749282495974337219   
 0.222522042695637753     0.9749276591186449509  
0.6234903946185663099   0.7818310097574695578   
 0.9009690323808488999   0.433883397574413817   /  

   \endpicture
  \hss \egroup
\endinsert

\bigskip
\bigskip
\noindent 3.4 {\bf {\sl The algebra $C(\Omega) \rtimes\Gamma$}}.

  As before, $\triangle$ is a triangle building,
$\Gamma$ is a group which acts simply transitively on the vertices of $\triangle$ and we
 identify $\Gamma$ with the vertices of its Cayley graph i.e. with the vertices of
$\triangle$.

We first partition the boundary $\Omega$ of $\triangle$.  Given generators $a, b \in P$
with $b \in \lambda(a)$ let $\Omega(a^{-1}, b)$ denote the set of elements $\omega \in
\Omega$ whose representative sector $S(\omega)$ has base chamber $\{ e, a^{-1}, b\}$.

\midinsert
  \line \bgroup \hss
\beginpicture
\setcoordinatesystem units <0.5cm, 0.866cm>
\put {$ Figure \ 3.4.1 $} [t] at 0 -0.6
\put {$\bullet$} at 0 0
\multiput {$\bullet$} at -1 1 *1 2 0 /
\put {$e$} [t] at 0 -0.2
\put {$a^{-1}$} [r,t] at -1.2 1
\put {$b$} [l,t] at 1.2 1
\putrule from -1 1 to 1 1
\setlinear \plot -1 1 0 0 1 1 /
\multiput {.} at 1 1 *40 .05 .05 /
\multiput {.} at -1 1 *40  -.05 .05  /
\put {$S(\omega)$}  at 0 2.5
\endpicture

  \hss \egroup
\endinsert
\bigskip

\noindent In the notation of $\S$3.2, $\Omega(a^{-1}, b) = \Omega^e (a^{-1}) \cap
\Omega^e(b)$. Since there is a unique sector based at $e$ representing
a point  $\omega \in \Omega$ , it follows that the sets $\Omega(a^{-1}, b)$
 are pairwise disjoint.  To find the number of possible $\Omega(a^{-1}, b)$ note that there are
$q^2+q+1$ choices for $a^{-1}$ and then $q+1$ choices for $b$, since each line contains
$q+1$ points.  It follows that $\Omega$ is partitioned into $(q+1)(q^2+q+1)$ open and
closed sets of the form $\Omega (a^{-1}, b)$.
The following key technical result will be used later.

\noindent {\bf Lemma 3.1}  Given generators $a, b \in P$
with $b \in \lambda(a)$  and any basic open set $\Omega^e(v)$ where
$v$ is a vertex of $\triangle$, there exists $k \in \G$ such that
 $k\Omega(a^{-1}, b) \subset \Omega^e(v)$

\noindent {\bf Proof}.   Assume that $v$ has left normal form
 $v =  x_1^{-1}x_2^{-1} \dots x_n^{-1}y_1y_2 \dots y_m$
, where $m \ge 1$. ( If $m = 0$ the argument is simpler.)
Since $b \in \lambda(a)$, we have $abc = 1$ for a (unique)
$c \in P$.

Consider the projective plane of neighbours of $e$.
In it there are $q+1$ lines containing $b$. Therefore there are
$q+1$ possible $z_1 \in P$ such that $b \in \lambda(z_1)$. Similarly, there are $q+1$
points on the line $\lambda(y_m)$; i.e. there are $q+1$ possible $z_2 \in P$ such that
$z_2 \in \lambda(y_m)$. Since $|P| = q^2 + q + 1 > 2(q + 1)$, we can choose  $z \in P$ such
that  $b \notin \lambda(z)$ and  $z \notin \lambda(y_m)$. It follows that the word
$ x_1^{-1}x_2^{-1} \dots x_n^{-1}y_1y_2 \dots y_mzb$ has no obvious cancellations and so, by
[9: Lemma 6.2], $|vzb| = |v| + 2$. Also, by the remark following [9: Lemma 6.2],
$|vza^{-1}| = |vzbc| \ge |v| +2$, and so $|vza^{-1}| = |v| + 2$.

Now let $\o_0 \in \O(a^{-1}, b)$.
Then the sector $S^{vz}(vz\omega_0)$ is  a subsector of a sector with base vertex $e$
which contains $v$. (See Figure 3.4.2.)

\midinsert
  \line \bgroup \hss
\beginpicture
\setcoordinatesystem units <0.5cm, 0.866cm>
\put {$Figure \  3.4.2 $} at 0 -1
\put {$\bullet$} at 0 0
\put {$\bullet$} at -1 2
\put {$\bullet$} at 0 3
\multiput {$\bullet$} at -1 4 *1 2 0 /
\put {$S^{vz}(vz\omega_0)$}  at 0 5
\put {$e$} [t] at 0 -0.1
\put {$v$} [l] at -0.8 1.9
\put {$vz$} [l] at 0.2 3
\put {$vzb$} [l] at 1.2 4
\put {$vza^{-1}$} [r] at -1.2 4
\setlinear \plot -4 4  0 0  4 4 /
\setlinear \plot -1 2  2 5 /
\setlinear \plot  0 3  -2 5 /
\multiput {.} at -1 4 *40  .05 0  /
\endpicture
  \hss \egroup
\endinsert

\bigskip

The sector with base vertex $e$ is parallel to  $S^{vz}(vz\omega_0)$  and so , by uniqueness, must be
 $S^{e}(vz\omega_0)$. Let $k = vz$. Then  $k\omega_0 \in \Omega^e(v)$.
Thus  $k\Omega(a^{-1}, b) \subset \Omega^e(v)$, as required.

If $m = 0$ then the argument is similar. We need only choose $z \in P$ with $b \notin \lambda(z)$
and $z \ne x_n$.
\ep
\bigskip
\bigskip

There is a natural well--defined action of $\Gamma$ on $\Omega$ coming from left
multiplication of every vertex in a sector by an element $g \in \Gamma$ [7 : $\S$2].
This induces an action $\alpha$ of $\Gamma$ by automorphisms of $C(\Omega)$ via
$$
\alpha(g) \ f(\omega) = f(g^{-1} \omega)
$$
\noindent By analogy with the tree case of $\S$2 we now define partial isometries
$s_{a^{-1},b}^+$, $s_{a^{-1},b}^-$ in $C(\Omega) \rtimes \Gamma$ which generate
$C(\Omega) \rtimes \Gamma$ and satisfy Cuntz--Krieger type relations.

First of all note that the characteristic function $1_{\Omega (a^{-1}, b)}$ is continuous on
$\Omega$ and so defines a projection $p_{a^{-1}, b}$ in $C(\Omega) \rtimes \Gamma$.

Define
$$
s_{a^{-1},b}^+ = b \sum  p_{c^{-1}, d} \eqno(1)
$$
where the sum is over those $(c, d)$ for which there is a sector of the form shown in Figure 3.4.3.
\bigskip
\bigskip
\midinsert
  \line \bgroup \hss
\beginpicture
\setcoordinatesystem units <1cm, 1.732cm>
\put {$Figure \ 3.4.3$} at 0 -1.5
\put {$a$} [r ] at -0.6 -0.5
\put {$b$} [ l] at 0.6 -0.5
\put {$x$} [ b] at 0 0.1
\put {$c$} [r ] at 0.4  0.5
\put {$d$} [ l] at 1.6  0.5
\multiput {\beginpicture
\setcoordinatesystem units <1cm, 1.732cm>
\arrow <10pt> [.2, .67] from  0.2 1 to 0 1
\arrow <10pt> [.2, .67] from  -0.7 0.7 to -0.5 0.5
\arrow <10pt> [.2, .67] from  0.3 0.3 to 0.5 0.5
\putrule from -1 1 to 1 1
\setlinear \plot -1 1 0 0 1 1 /
\endpicture}  at   0 0  1 1 /
\setlinear \plot -1 0  -2.5 1.5 /
\setlinear \plot  2 1  2.5 1.5  /
\endpicture
  \hss \egroup
\endinsert

\bigskip
\noindent The possible $(c, d)$ are determined as follows.  Firstly $d \not\in \lambda (b)$
and so there are  $(q^2+q+1)$ \hfill \break
$- (q+1) = q^2$ choices for $d$.  Once $d$ is chosen, we have $x \in \lambda(c)$
and $d \in \lambda (c)$, which determines $c$ uniquely.  Therefore there are
$q^2$ terms in the sum (1).  Denote by $A_{a^{-1}, b}^+$ the set of such $(c,
d)$.
\noindent Now $ s_{a^{-1}, b}^+$   is a partial isometry with
 $$ s_{a^{-1}, b}^+ \  s_{a^{-1}, b}^{+*}  =  p_{a^{-1}, b}$$

and

$$s_{a^{-1}, b}^{+*} \  s_{a^{-1}, b}^+ = \sum_{(c, d) \in A_{a^{-1},b}^+}
s_{c^{-1}, d}^+ \   s_{c^{-1}, d}^{+*} .\eqno(2)
$$

\noindent These partial isometries therefore generate a Cuntz--Krieger algebra [10].

The action of $s_{a^{-1},b}^+$ may be represented pictorially as in Figure 3.4.4.
\bigskip \bigskip

\midinsert
 \line \bgroup \hss
\beginpicture
\setcoordinatesystem units <0.5cm,0.866cm>    
\put {$Figure\ 3.4.4$} at  8  -0.5
\put {$\bullet$} at  0 0
\multiput {$\bullet$} at  -1 1  *1  2 0 /
\put {$\longrightarrow$} at 8 2
\put {$b$}   [b]  at  8 2.2
\put {$c^{-1}$} [r] at  -1.2 1
\put {$d$}      [l] at  1.2 1
\put {$e$}      [t] at   0 -0.2
\putrule from -1 1 to 1 1
\setlinear \plot -3.5 3.5  0 0  3.5 3.5 /
\setcoordinatesystem units <0.5cm,0.866cm>  point at -15 0 
\put {$\bullet$} at  1 1
\multiput {$\bullet$} at  0 2  *1  2 0 /
\put {$a^{-1}$} [r] at  -1.2 1
\put {$b$}      [l] at  1.2 1
\put {$bc^{-1}$} [r] at -0.2 2
\put {$bd$}      [l] at  2.2 2
\put {$e$}       [t] at  0 -0.2
\setlinear \plot  -1.5 3.5  1 1  3.5 3.5 /
\putrule from  0 2  to 2 2
\setdashes
\putrule from -1 1 to 1 1
\setdashes
\setlinear \plot  -1 1    0 2  /
\setdashes
\setlinear \plot -3.5 3.5  0 0  1 1 /
\endpicture
 \hss \egroup
\endinsert

\smallskip

\noindent Note that the unbroken sector on the right with base vertex $b$ is parallel to the whole sector
 and hence
is a representative of an element in $\Omega (a^{-1}, b)$.

The following relations are also clearly satisfied
$$
\sum_{(c,d)\in A^+_{a^{-1},b}} b\  p_{c^{-1}, d}\  b^{-1} = p_{a^{-1}, b} .
\eqno(3) $$
\noindent Similarly we define
$$
s_{a^{-1}, b}^- = a^{-1} \sum p_{c^{-1}, d}, \eqno(4)
$$
where the sum is over the set $A_{a^{-1}, b}^-$ of those  $(c, d)$ for which there is a sector
as inFigure 3.4.5
\bigskip
\bigskip

\midinsert
  \line \bgroup \hss
\beginpicture
\setcoordinatesystem units <1cm, 1.732cm>
\put {$Figure\ 3.4.5$} at  0  -1.5
 \put {$a$} [r ] at -0.6 -0.5
\put {$b$} [ l] at 0.6 -0.5
\put {$c$} [r ] at -1.6  0.5
\put {$d$} [ l] at -0.4  0.5
\multiput {\beginpicture
\setcoordinatesystem units <1cm, 1.732cm>
\arrow <10pt> [.2, .67] from  0.2 1 to 0 1
\arrow <10pt> [.2, .67] from  -0.7 0.7 to -0.5 0.5
\arrow <10pt> [.2, .67] from  0.3 0.3 to 0.5 0.5
\putrule from -1 1 to 1 1
\setlinear \plot -1 1 0 0 1 1 /
\endpicture}  at   0 0  -1 1 /
\setlinear \plot -2 1  -2.5 1.5 /
\setlinear \plot  1 0  2.5 1.5  /
\endpicture
  \hss \egroup
\endinsert

\bigskip
\bigskip
where $a \not\in \lambda(c)$.  Again there are $q^2$ terms in the sum.

\noindent Then $s_{a^{-1}, b}^-$ is a partial isometry with
$$
s_{a^{-1}, b}^- \ \ s_{a^{-1}, b}^{-^{*}} = p_{a^{-1}, b}
$$
and
$$
 s_{a^{-1}, b}^{-^{*}} \ \ s_{a^{-1}, b}^- = \sum_{(c,d)\in A_{a^{-1},b}^-} \
s_{c^{-1}, d}^- \ s_{c^{-1}, d}^{-^{*}} \eqno(5)
$$
so that these partial isometries also generate a Cuntz--Krieger algebra.  There are also
covariance relations analogous to (3).

The action  of $s_{a^{-1}, b}^-$ can be represented by Figure 3.4.6.
\bigskip \bigskip

\midinsert
 \line \bgroup \hss
\beginpicture
\setcoordinatesystem units <0.5cm,0.866cm>    
\put {$Figure \ 3.4.6$} at  8  -0.5
\put {$\bullet$} at  0 0
\multiput {$\bullet$} at  -1 1  *1  2 0 /
\put {$\longrightarrow$} at 8 2
\put {$a^{-1}$}   [b]  at  8 2.2
\put {$c^{-1}$} [r] at  -1.2 1
\put {$d$}      [l] at  1.2 1
\put {$e$}      [t] at   0 -0.2
\putrule from -1 1 to 1 1
\setlinear \plot -3.5 3.5  0 0  3.5 3.5 /
\setcoordinatesystem units <0.5cm,0.866cm>  point at -15 0 
\put {$\bullet$} at  -1 1
\multiput {$\bullet$} at  -2 2  *1  2 0 /
\put {$a^{-1}$} [r] at  -1.2 1
\put {$b$}      [l] at  1.2 1
\put {$a^{-1}c^{-1}$} [r] at -2.2 2
\put {$a^{-1}d$}      [l] at  0.2 2
\put {$e$}       [t] at  0 -0.2
\setlinear \plot  -3.5 3.5  -1 1  1.5 3.5 /
\putrule from  -2 2  to 0 2
\setdashes
\putrule from -1 1 to 1 1
\setdashes
\setlinear \plot  -1 1    0 0  3.5 3.5  /
\setdashes
\setlinear \plot  2 2  1 1  0 2 /
\endpicture
 \hss \egroup
\endinsert

\smallskip
Denote by $C^*(s^+)$, $C^*(s^-)$ the $C^*$--subalgebras of $C(\Omega) \rtimes \Gamma$
generated by the operators $s_{a^{-1}, b}^+$, $s_{a^{-1}, b}^-$ respectively.  Clearly
$$
C^*(s^+) \cap C^*(s^-) \supset \left \{p_{a^{-1}, b} : b \in \lambda(a) \right \}.
$$

We are interested in the $C^*$--algebra $C^*(s^\pm)$ which is generated by both sets of
operators.

\proclaim  Theorem. $C^*(s^\pm) = C(\Omega) \rtimes \Gamma$.

\noindent {\bf Proof}.  First of all we show that any generator $b \in \Gamma$ lies in
$C^*(s^\pm)$.  We do this by demonstrating the following explicit formula for $b$ in terms
of elements of $C^*(s^\pm)$.
$$
b = \sum_a  s_{a^{-1}, b}^+ + \sum_k  s_{b^{-1}, k}^{-^{*}} + \sum  s_{t^{-1}, f}^-
s_{h^{-1}, s}^{+^{*}} . \eqno(6)
$$

We consider the three sums separately.

\noindent (A) {\bf The first sum}  is over those $a \in P$ such that $b \in \lambda(a)$.  Thus
the sum contains $q+1$ terms and the initial projection of each term is the sum of the
projections of the form $p_{c^{-1}, d} = 1_{\Omega(c^{-1}, d)}$ satisfying the additional
 requirement that $(c, d) \in A_{a^{-1}, b}^+$.
[See equation (2).]   As previously noted, there are $q^2$ such projections, and so the initial
space of the first sum is the sum of $q^2(q+1)$ projections of the form $p_{c^{-1}, d}$.

\noindent (B){\bf The second sum} is over $k \in \lambda(b)$.  There are $q+1$ such k's and
the initial projection of $s_{b^{-1}, k}^{- *}$ is $p_{b^{-1}, k}$.  Thus the initial space of the
second sum is the sum of $q+1$ projections of the form $p_{b^{-1}, k}$.
  Note that $p_{b^{-1}, k}$ is different from the projections $p_{c^{-1}, d}$ in (A), because
$k \in \lambda(b)$ whereas $d \notin \lambda(b)$; so that  $k \ne d$.

\noindent (C) {\bf The third sum}  is over those $t, f, h, s$ for which the following diagram is
possible.
\bigskip


\midinsert
  \line \bgroup \hss
\beginpicture
\setcoordinatesystem units <1cm, 1.732cm>
\put {$Figure 3.4.7$} at  1  -1.6
\put {$h$} [r ] at -0.6 -0.5
\put {$s$} [ l] at 0.5 -0.6
\put {$t$} [r ] at 1.5  -0.6
\put {$f$} [ l] at 2.6  -0.5
\put {$b$} [t]  at 1 -1.2
\arrow <10pt> [.2, .67] from 1.2 -1 to 1 -1
\putrule from  0 -1  to  2 -1
\multiput {\beginpicture
\setcoordinatesystem units <1cm, 1.732cm>
\arrow <10pt> [.2, .67] from  0.2 1 to 0 1
\arrow <10pt> [.2, .67] from  -0.7 0.7 to -0.5 0.5
\arrow <10pt> [.2, .67] from  0.3 0.3 to 0.5 0.5
\putrule from -1 1 to 1 1
\setlinear \plot -1 1 0 0 1 1 /
\endpicture}  at   0 0  2 0 /
\endpicture
  \hss \egroup
\endinsert

\medskip
\noindent That is
$$
\eqalign {
& s \in \lambda(b)   \quad \quad \quad \quad (q+1 \hbox { choices for } s) \cr
&b = t^{-1}s^{-1} \cr
&h \not= b \hbox { and } s \in \lambda (h) \quad \quad (q \hbox { choices for } h) \cr
&f \not= b  \hbox { and } f \in \lambda (t) \quad \quad (q \hbox { choices for } f) \cr}
$$

\noindent Therefore the third sum has $q^2(q+1)$ terms and yields a partial isometry whose
initial projection equals the sum of $q(q+1)$ projections  $p_{h^{-1}, s}$ .

  Now  $p_{h^{-1}, s}$ differs from any projection $p_{c^{-1}, d}$ of (A) since $ s \in \lambda (b)$
 in the first case but $ d \notin \lambda (b)$ in the second; and it differs from any $p_{b^{-1}, k}$
 in (B) since $h \ne b$. The initial
projections involved in (A), (B), and (C) are therefore all different.

For all three sums the total number of initial projections of the form $p_{x^{-1}, y}$ is
$$
q^2(q+1) + (q+1) + q(q+1) = (q+1)(q^2+q+1).
$$

\noindent Thus the initial projection is the identity operator. Similarly the final projection is the sum of $(q+1) +
q^2(q+1) + q(q+1) = (q+1)(q^2+q+1)$ projections which sum to the identity.

On each initial projection the sum on the right of (6) is multiplication by $b$ and hence (6) is
established.  Since each generator of $\Gamma$ lies in $C^*(s^\pm)$, the same is true for
any element of $\Gamma$.

We can now complete the proof.  Recall that
$$
1_{\Omega(a^{-1}, b)} = p_{a^{-1}, b} = s_{a^{-1}, b}^+ \  s_{a^{-1}, b}^{+ *}
$$
and so $1_{\Omega(a^{-1}, b)}$ lies in $C^*(s^\pm)$.
\noindent Therefore $1_{g\Omega(a^{-1}, b)} = \alpha(g) \ 1_{\Omega(a^{-1}, b)} = g
1_{\Omega(a^{-1}, b)} \ g^{-1}$ lies in $C^*(s^\pm)$ for each $g \in \Gamma$, since
$\Gamma \subset C^*(s^\pm)$.  Similarly $1_{g\left (\Omega \setminus \Omega(a^{-1},
b)\right)} \in C^*(s^\pm), g \in \Gamma$.  We claim that the family of $g
\Omega (a^{-1}, b)$ and $g\left ( \Omega \setminus \Omega (a^{-1}, b) \right )$ forms a
base for the topology of $\Omega$.  The Stone--Weierstrass theorem
then shows that $C(\Omega) \subset
C^*(s^\pm)$.  Therefore $C(\Omega) \rtimes \Gamma = C^*(s^\pm)$.

In order to prove our claim,
let $\Omega^e(v)$ be any basic open set, and let
 $\omega\in\Omega^e(v)$.  Define $a$ and $b$
by $S^v(\omega)\in\Omega^v(va^{-1})\cap\Omega^v
(vb)$. In other words, $a$ and $b$ are defined
by the labels on the left and right edges of
the base triangle in $S^v(\omega)$.
Since $S^v(\omega)\subseteq S^e(\omega)$,
we have $|va^{-1}|=|vb|=|v|+1$.  Then, as in the
proof of Lemma 3.1, $v\Omega(a^{-1},b)\subseteq
\Omega^e(v)$, and by construction, $\omega\in
v\Omega(a^{-1},b)$.
 \ep

\noindent {\bf Remark}. Since $C^*(s^+) \cap C^*(s^-) \supset \left \{p_{a^{-1},b} :
 b \in \lambda(a) \right \}$ , the above proof shows that $C(\Omega) \subset
\alpha(\Gamma)(C^*(s^+) \cap C^*(s^-))$.
\bigskip
\bigskip

\noindent 3.5 {\bf {\sl Weak Commutativity}}.

  A natural project is to study the algebras  $C^*(s^\pm)$ by analogy with the
Cuntz-Krieger algebras. One would hope, for example that the relations  (2)
and  (5)  defining the Cuntz-Krieger algebras  $C^*(s^+)$ and  $C^*(s^-)$
together with some connecting relations would characterize the generated  $C^*-$algebra
up to canonical isomorphism.We show a weak form of commutativity between $C^*(s^+)$ and
$C^*(s^-)$.

Suppose that  $s^+_{a^{-1},b}s^-_{c^{-1},d} \ne 0$, that is  $(c,d) \in A^+_{a^{-1},b}$.
Then
$$s^+_{a^{-1},b}s^-_{c^{-1},d} = bc^{-1}\sum_{(e,f) \in A^-_{c^{-1},d}} p_{e^{-1},f}$$

Similarly, if  $s^-_{a^{-1},b}s^+_{g^{-1},h} \ne 0$, that is $(g,h) \in A^-_{a^{-1},b}$,
then
$$s^-_{a^{-1},b}s^+_{g^{-1},h} = a^{-1}h \sum_{(e,f) \in A^+_{g^{-1},h}} p_{e^{-1},f}$$

\midinsert
  \line \bgroup \hss
\beginpicture
\setcoordinatesystem units <1cm, 1.732cm>
\setplotarea  x from -3 to 3,  y from -2 to 2
\put {$Figure \  3.5.1$} at  0 -1.5
\put {$a$} [r ] at -0.6 -0.5
\put {$b$} [ l] at 0.6 -0.5
\put {$c$} [r ] at 0.4  0.5
\put {$d$} [ l] at 1.6  0.5
\put {$g$} [r ] at -1.5  0.4
\put {$h$} [ l] at -0.4  0.5
\put {$e$} [r ] at -0.6 1.5
\put {$f$} [ l] at 0.6 1.5
\multiput {\beginpicture
\setcoordinatesystem units <1cm, 1.732cm>
\arrow <10pt> [.2, .67] from  0.2 1 to 0 1
\arrow <10pt> [.2, .67] from  -0.7 0.7 to -0.5 0.5
\arrow <10pt> [.2, .67] from  0.3 0.3 to 0.5 0.5
\putrule from -1 1 to 1 1
\setlinear \plot -1 1 0 0 1 1 /
\endpicture}  at   0 0  1 1  -1 1  0 2 /
\setlinear \plot -2 1  -3.5 2.5 /
\setlinear \plot  2 1  3.5 2.5  /
\setlinear \plot -1 2  -1.5 2.5 /
\setlinear \plot  1 2   1.5 2.5 /
\endpicture
  \hss \egroup
\endinsert

Now suppose that $bc^{-1} = a^{-1}h$ and that $(e,f) \in A^-_{c^{-1},d} \cap A^+_{g^{-1},h}$
[there are  $q$ possible such $(e,f)$].
This means that there is a diagram as in Figure 3.5.1

Then $s^+_{a^{-1},b}s^-_{c^{-1},d} p_{e^{-1},f} = s^-_{a^{-1},b}s^+_{g^{-1},h} p_{e^{-1},f}$

Summing over all possibilities gives
$$\sum_ds^+_{a^{-1},b}s^-_{c^{-1},d} = \sum_gs^-_{a^{-1},b}s^+_{g^{-1},h}.  \eqno(7)$$
For fixed $a$,$b$,$c$,$h$, there are $q$ terms in each sum.
\bigskip
\bigskip
\bigskip
\noindent 4. {\bf The action of $\Gamma$ on $\O$}.
\bigskip
\noindent 4.1 {\bf Minimality of the action.}

By definition, the action of $\G$ is {\it minimal} if, for all $\o \in \O$, $\G\o$ is
dense in $\O$.
\bigskip

\noindent {\bf Proposition 4.1.1} The action of $\Gamma$ on $\Omega$ is minimal.

\noindent {\bf Proof}. Let $\omega_0 \in \Omega$ and let $\Omega^e(v) = \{\omega \in \Omega:
v \in S^e(\omega)\}$ be a basic open set in $\Omega$, where $v$ is a vertex of $\Delta$
(i.e. $v \in \Gamma$)  and $S^e(\omega)$ denotes the unique sector in the class $\omega$
having base vertex $e$.
We must find $k \in \Gamma$ such that $k\omega_0 \in \Omega^e(v)$.
 Let the sector $S^e(\omega_0)$
have base chamber $\{e,a^{-1},b\}$, where $b \in \lambda(a)$; that is $\o_0 \in \O(a^{-1}, b)$.
By Lemma 3.1, there exists $k \in \G$ such that  $k\Omega(a^{-1}, b) \subset \Omega^e(v)$.
In particular $k\omega_0 \in \Omega^e(v)$.
\ep

\bigskip
\bigskip

\bigskip
\noindent 4.2. {\bf The action of $\Gamma$ on $\Omega$ is amenable}.

By analogy with the Cuntz-Krieger algebras one would hope that the algebras
$C^*(s^{\pm})$ are nuclear. We approach this problem by using the results
of [2]. If we can show that the action of $\Gamma$ on $\Omega$ is amenable
in the sense of [2] then  $C^*(s^\pm) = \cross$  is isomorphic to the reduced crossed
product $\crossr$ [2: Proposition 4.8]  and is nuclear [2: Th\'eor\`eme 4.5].

We verify that the action of $\G$ on $\O$ is amenable by constructing a sequence $f_i$
of continuous real valued functions with compact support on $\G \times \O$ such that

(a)   \quad  $\sum_{t \in \G} |f_i(t,\o)|^2 = 1$ \quad   for all $\o \in \O$ and $i \in \NN$

(b)   \quad  $lim_i \sum_{t \in \G}\overline{f_i(t,\o)}f_i(s^{-1}t,s^{-1}\o) = 1$ \quad
uniformly on $\O$ for each $s \in \G$

(See Th\'eor\`eme 4.9(h) and Remarque 4.10 in [2].)

\bigskip

 Denote by $\o_{m,n}$ the vertices of the sector $S^e(\o)$.
\bigskip
\midinsert                        
 \line \bgroup \hss
\beginpicture
\setcoordinatesystem units <0.5cm,0.866cm>   
\put {$\bullet$} at 0 0
\multiput {$\bullet$} at -1 1 *1 2 0 /
\multiput {$\bullet$} at -2 2 *2 2 0 /
\put {$Figure \  4.2.1$} at 0 -1
\put{$e$}     [t]   at  0  -0.2
\put{$\omega_{1,1}$} [b r] at  -0.2 2.05
\put{$\omega_{0,1}$} [r]  at -1.2 1
\put{$\omega_{0,2}$}  [r]  at -2.2 2
\put{$\omega_{1,0}$}  [l] at  1.2 1
\put{$\omega_{2,0}$}  [l] at 2.2 2
\putrule from -1 1 to 1 1
\putrule from -2 2 to 2 2
\setlinear \plot -2 2  0 0  2 2  /
\setlinear \plot -1 1  0 2  1 1 /
\setdashes
\setlinear \plot -3.5 3.5  -2 2  -1 3  0 2  1 3  2 2  3.5 3.5 /
\putrule from  -3 3  to  3 3
\endpicture
  \hss \egroup
\endinsert

 For $\o \in \O$, define  $f_i(\cdot ,\o)$ to be the characteristic function of the triangular
region
$\{ t \in S^e(\o) : |t| \le i-1\}$, normalized so that  (a)  holds. i.e.
 $\|f_i(\cdot ,\o)\|_2 = 1$.

By definition,
$$
f_i(t,\o) = \cases {(i(i+1)/2)^{-1/2}  \, & if  $t \in S^e(\o)$ and $ |t| \le i - 1$ \cr
0 &  otherwise. \cr}
$$

and there are exactly $i(i+1)/2$ elements $t \in \G$ for which $f_i(t,\o) \ne 0$. Also,
since $s^{-1}t \in S^e(s^{-1}\o)$ if and only if $t \in S^s(\o)$, we have
$$
f_i(s^{-1}t,s^{-1}\o) = \cases {(i(i+1)/2)^{-1/2}  \, & if  $t \in S^s(\o)$ and
 $ |s^{-1}t| \le i - 1$ \cr
0 &  otherwise. \cr}
$$

We estimate the number $N_i$ of vertices  $t \in S^e(\o) \cap S^s(\o)$ such that
$|t| \le i-1$ and $|s^{-1}t| \le i-1$.  Assume that  $i \ge |s|$  and note that by
[7: Corollary 2.3], $\o_{m,n} \in S^e(\o) \cap S^s(\o)$ for $m,n \ge |s|$.

\midinsert
  \line \bgroup \hss
\beginpicture
\setcoordinatesystem units <0.5cm, 0.866cm>
\setplotarea  x from -7 to 8,  y from -1 to 9
\put {$Figure \ 4.2.2$} at 0 -1
\put {$\bullet$} at 0 0
\put {$\bullet$} at 1 0.2
\put {$e$} [t] at 0 -0.2
\put {$s$} [t] at 1  0
\put {$\bullet$} at 0 2.4
\put {$\bullet$} at 2.8 5.2
\put {$\bullet$} at 0.4 7.6
\put {$\bullet$} at -2.4 4.8

\put {$v$} [b] at 0 2.6
\put {$w$} [r] at 2.6 5.2
\put {$x$} [t] at 0.4 7.4
\put {$y$} [l] at -2.2 4.8

\put {$\omega_{0,i}$} [r] at -4.2 4
\put {$\omega_{i,0}$} [r] at 3.8 4
\put {$\omega_{0,|s|}$} [r] at -1.4 1.2
\setlinear \plot -6 6  0 0  6 6  /
\setlinear \plot -5 6.2  1 0.2  7 6.2  /
\setlinear \plot -4 4  0 8  4 4  /
\setlinear \plot -3 4.2  1 8.2  5 4.2  /
\setdashes
\setlinear \plot 1.2 1.2  -3.8 6.2  /  
\setlinear \plot -1.2 1.2  3.8 6.2  /  
\endpicture
  \hss \egroup
\endinsert
\bigskip

Let $v = \o_{|s|,|s|}$ and $w= \o_{i,|s|}$. In Figure 4.2.2, the convex hull of $\{x,y,v,w\}$
 consists of vertices which lie in $S^e(\o) \cap S^s(\o)$.

Now $|\o_{m,n}| \le i-1$ whenever $m+n \le i-1$, and $|s^{-1}\o_{m,n}| \le i-1$
 whenever $m+n \le i- |s| -1$ , so clearly
$ N_i \ge (i - 3|s|)(i - 3|s| + 1)/2$.
Therefore
$\overline{f_i(t,\o)}f_i(s^{-1}t,s^{-1}\o) = (i(i+1)/2)^{-1}$
 for at least  $(i - 3|s|)(i - 3|s| + 1)/2$ values of $t$.
It follows that
$$
1 \ge \sum_{t\in\G}\overline{f_i(t,\o)}f_i(s^{-1}t,s^{-1}\o)
\ge (i(i+1)/2)^{-1}(i - 3|s|)(i - 3|s| + 1)/2 \rightarrow 1$$
as $i \rightarrow \infty$.
Condition (b) above is then satisfied and therefore the action of $\G$ on $\O$ is
amenable. We have therefore proved

\noindent {\bf Proposition 4.2.1}. The algebra $C^*(s^\pm)$ is nuclear.
\vfill\eject

\noindent 4.3. {\bf The action of $\G$ on $\O$ is topologically free}

  The proof of the stated assertion requires a preliminary result , which is of independent
interest. We refer to the appendix for further information on the concepts used in
its proof.

\noindent {\bf Proposition 4.3.1} Let $V$ be an open subset of $\O$. There exists an
apartment each of whose six boundary points lies in $V$.

\noindent {\bf Proof}.
Identify the boundary $\Omega$  with the set of
$\pi/3$-angled sectors based at $e$. That is, identify $\o\in\O$ with $S^e(\o)$.
Any open set $V\subseteq\Omega$ contains a set of the form
$$
\Omega(T)=\{ \omega \in \Omega : T \subset S^e(\o) \}
$$
for some (big) equilateral trangle $T$ based at $e$. 

We shall find an apartment each of whose 6 boundary points lies in
$\Omega(T)$.  To start, for reasons that become clear only later,
we add one layer to $T$, choosing it arbitrarily.  Then we add a
big ``upside down'' equilateral triangle to the extended $T$, calling
the whole resulting diamond shaped figure $T'$. See Figure 4.3.1.
The tip vertex of $T'$ we label $v_1$.

\midinsert
  \line \bgroup \hss
\beginpicture
\setcoordinatesystem units <0.5cm, 0.866cm>  
\setplotarea  x from -7 to 8,  y from -1 to 9
\put {$Figure \ 4.3.1$} at 0 -1
\put {$\bullet$} at 0 0
\put {$\bullet$} at 0 8
\put {$e$} [t] at 0 -0.2
\put {$v_1$} [b] at 0 8.2
\put {$T$} at  0  2.5
\putrule from -3.5 3.5  to  3.5 3.5
\putrule from -4 4 to 4 4
\setlinear \plot -4 4  0 0  4  4  /
\setlinear \plot -4 4  0 8  4 4    /
\endpicture
  \hss \egroup
\endinsert
\bigskip

By an element of $\triangle$ we mean any vertex, edge, or triangle
(chamber) of $\triangle$.  For any such element $E$ there exist sectors
 based at $e$ which contain $E$.  Moreover, although there are
uncountably many such sectors, in each of them $E$ lies in the same
relative location.  This gives a retraction, $r_e$, of the entire
building onto an abstract sector.  By abstract, we mean a sector
 which is specially intended as the image of $r_e$, instead of
being some particular sector of the building.  To say that $E$
``lies'' in a given position is to say that $r_e(E)$ takes some
particular value.  The particular value is easily indicated on a
figure.

 Each edge is contained in exactly
$q+1$ triangles.  If we retract an edge, which does not lie on a ray
through $e$, and its $q+1$ containing
triangles, then one of the triangles will lie below the edge and the
other $q$ will lie above it.  This applies whether the edge lies
horizontally or at a slant, and it is so because there is a unique
stretched gallery of a given type from $e$ to the edge.  The last
triangle of that stretched gallery is the one which lies below the
edge.

Let $A$ be a triangle emanating from $v_1$ which lies as indicated
in Figure 4.3.2.
\bigskip
\bigskip
\midinsert
  \line \bgroup \hss
\beginpicture
\setcoordinatesystem units <0.5cm, 0.866cm>
\setplotarea  x from -7 to 8,  y from -1 to 10
\put {$Figure \  4.3.2$} at 0 -1
\put {$\bullet$} at 0 0
\put {$\bullet$} at 0 8
\put {$e$} [t] at 0 -0.2
\put {$v_1$} [r] at -0.2 8
\put {$T'$} at  0  4
\put {$A$} at  0 8.7
\put {$W_1$} at  0  10
\putrule from -1 9 to 1 9
\setlinear \plot 1 9  0 8  -1 9  /
\setlinear \plot -4 4  0 0  4  4  0 8  -4 4 /
\setdashes
\setlinear \plot -6 6 -4 4  /
\setlinear \plot 4 4  6 6  /
\setlinear \plot   1 9  2.5 10.5  /
\setlinear \plot  -1 9  -2.5 10.5  /
\endpicture
  \hss \egroup
\endinsert

  Any sector $W_1$ based at $v_1$ with first
triangle $A$ will lie as indicated in Figure 4.3.2 .
(Proof: one may construct a gallery of stretched type from $e$ to
any element $E$ of $W_1$ out of a stretched gallery from $e$ to $A$
and another from $A$ to $E$.  Any  gallery of stretched type (i.e. having
the same type as a stretched gallery) is in
fact stretched, and therefore lies in some apartment. This
establishes the value of $r_e(E)$, that is, it establishes
where $E$ lies.)

One can fill in the figure consisting of $T'$ and $W_1$ to a unique
sector based at $e$.  To see this, just  fill in one layer
at a time starting from the boundary of the given figure and work towards the
boundary of the desired figure.
For example, starting at the vertex $v_1$, the chambers $A$ and $T_1$ in
Figure 4.3.2(a) are determined. Each vertex of the building which is adjacent
to $v_1$ corresponds to a point or line of a finite projective plane.
It follows that $z$ is uniquely determined, since it is incident with both
$x$ and $y$. Similarly $w$ is uniquely determined. All the chambers in
Figure 4.3.2(a) are now determined. Continue in this way, working next
from the the vertex  $x$ instead of $v_1$ .

\midinsert
  \line \bgroup \hss
\beginpicture
\setcoordinatesystem units <1cm, 1.732cm>
\setplotarea  x from -2.5 to 2.5,  y from -1.5 to 1.5
\put {$Figure \  4.3.2(a)$} at 0 -1.5
\put {$\bullet$} at 0 0
\put {$\bullet$} at 1 1
\put {$\bullet$} at 1 -1
\put {$\bullet$} at -1 1
\put {$\bullet$} at -1 -1
\put {$v_1$} [b,l] at 0.2 0.05
\put {$x$} [l] at  1.2 1
\put {$y$} [l] at  1.2 -1
\put {$z$} [l] at  2.2 0
\put {$w$} [r] at -2.2 0
\put {$A$} at  0 0.7
\put {$T_1$} at  0 -0.7
\putrule from -1 1 to 1 1
\putrule from -2 0 to 2 0
\putrule from -1 -1 to 1 -1
\setlinear \plot 0 0  1 1  2 0  1 -1  0 0  -1 1  -2 0  -1 -1  0 0 /
\endpicture
  \hss \egroup
\endinsert

It follows that the point of $\Omega$ represented by $W_1$ belongs to
$\Omega(T)$, indeed to
$$
\Omega(T')=\{ \omega\in\Omega :
                 T' \subset S^e(\omega) \}
$$

 Consider the {\it star} of the vertex $v_1$. That is the union of all the
 triangles of the building which contain $v_1$. In this star find a hexagon
which has four triangles  $A_1$ , \dots , $A_4$  which all lie as $A$ does.
  Choose $A_1$ from among
the $q^3$ possibilities.  There are then $q$ triangles which share
with $A_1$ the edge that would make them possible $A_2$'s.  Choose
one of the $q-1$ of these which lies as $A$ rather than as $D$ in
Figure 4.3.3.
\midinsert
  \line \bgroup \hss
\beginpicture
\setcoordinatesystem units <0.5cm, 0.866cm>
\setplotarea  x from -7 to 8,  y from -1 to 11
\put {$Figure \  4.3.3$} at 0 -1
\put {$\bullet$} at 0 0
\put {$\bullet$} at 0 8
\put {$e$} [t] at 0 -0.2
\put {$T'$} at  0  4
\put {$A$} at  0 8.7
\put {$B$} at -1 8.3
\put {$C$} at  -2 8.7
\put {$F$} at -3 8.3
\put {$D$} at  1 8.3
\putrule from -3 9 to 1 9
\putrule from  -4 8 to 2 8
\setlinear \plot 1 9  2 8  /
\setlinear \plot -4 4  0 0  4  4  0 8  -4 4  /
\setlinear \plot -4 8  -3 9  -2 8  -1 9  0 8  1 9  /
\setdashes
\setlinear \plot -6 6 -4 4  /
\setlinear \plot 4 4  6 6  /
\endpicture
  \hss \egroup
\endinsert
\bigskip

  Similarly, there are again $q-1$ choices for $A_3$, and
again for $A_4$.  The $A_1$ , \dots , $A_4$ so chosen form a stretched
gallery, so they lie together in some apartment; in particular they lie
in some hexagon (Figure 4.3.4).
\bigskip \bigskip
\midinsert
  \line \bgroup \hss
\beginpicture
\setcoordinatesystem units <1cm, 1.732cm>
\setplotarea  x from -2 to 2,  y from -1 to 1
\put {$Figure \  4.3.4$} at 0 -1.5
\put {$\bullet$} at 0 0
\put {$v_1$} [b,l] at 0.2 0.05
\put {$A_2$} at  1.1 0.4
\put {$A_1$} at  0 0.7
\put {$A_3$} at  1.1 -0.4
\put {$A_4$} at  0 -0.7
\put {$B_1$} at -1.1 0.4
\put {$B_2$} at -1.1 -0.4
\putrule from -1 1 to 1 1
\putrule from -2 0 to 2 0
\putrule from -1 -1 to 1 -1
\setlinear \plot 0 0  1 1  2 0  1 -1  0 0  -1 1  -2 0  -1 -1  0 0 /
\endpicture
  \hss \egroup
\endinsert

The remaining chambers $B_1$ , $B_2$  of the hexagon
 might lie as $A$ or they might
lie as $B$ .  It is not possible for one to lie one way
and the other to lie the other.  If they both lie as $A$, then we
can construct the desired apartment by expanding out from the hexagon.
For then each of the six sectors in this apartment based at $v_1$
lies as $W_1$ ( Figure 4.3.2 )
 and so represents a point of $\O$ belonging to $\O(T)$.
  So we assume that $B_1$ and $B_2$ both lie as $B$ .

Let us expand the hexagon in Figure 4.3.4 as follows to make an
 asymmetric hexagon of ten triangles
( Figure 4.3.5 ) which we will take as part of our apartment.  First
choose $C_1$ in any of $q$ ways.  All of these lie as $C$  in Figure 4.3.3 ,
 because $B_1$ lies below the edge it shares with the potential
$C_1$'s.  Next choose as $C_1'$ any of the $q-1$ triangles bordering
$C_1$ and lying as $C$ rather than as $F$.  The choices of $C_1$ and
$C_1'$ determine $C_2'$ and $C_2$ uniquely.  Just like $C_1$, $C_2$
necessarily lies as $C$.  This then implies that $C_2'$ also lies as
$C$, because we know where two of its edges lie.

\midinsert
  \line \bgroup \hss
\beginpicture
\setcoordinatesystem units <1cm, 1.732cm>
\setplotarea  x from -3 to 3,  y from -2 to 2
\put {$Figure \  4.3.5$} at 0 -2
\put {$\bullet$} at 0 0
\put {$\bullet$} at -2 0
\put {$v_1$} [b,l] at 0.2 0.05
\put {$v_2$} [b,l] at -1.8 0.05
\put {$A_2$} at  1.1 0.4
\put {$A_1$} at  0 0.7
\put {$A_3$} at  1.1 -0.4
\put {$A_4$} at  0 -0.7
\put {$B_1$} at -1.1 0.4
\put {$B_2$} at -1.1 -0.4
\put {$C_1$} at  -2 0.7
\put {$C_1'$} at -3.1 0.4
\put {$C_2$} at -2 -0.7
\put {$C_2'$} at -3.1 -0.4
\put {$W_2$} at -4.5 1
\put {$W_3$} at -4.5 -1
\putrule from -3 1 to 1 1
\putrule from -4 0 to 2 0
\putrule from -3 -1 to 1 -1
\setlinear \plot 0 0  1 1  2 0  1 -1  0 0  -1 1  -2 0  -1 -1  0 0 /
\setlinear \plot  -2 0  -3 -1  -4 0  -3 1  -2 0 /
\setdashes
\putrule from -6 0 to -4 0
\setlinear \plot -4 2  -3 1  /
\setlinear \plot -4 -2  -3 -1  /
\endpicture
  \hss \egroup
\endinsert

Extend the asymmetric hexagon to an apartment ${\cal A}$ , and
consider the sector $W_2$ in that apartment.  Let $v_2$ be the
base vertex of $W_2$. Then  $W_2$ necessarily lies as shown in Figure 4.3.6

\midinsert
  \line \bgroup \hss
\beginpicture
\setcoordinatesystem units <0.5cm, 0.866cm>
\setplotarea  x from -7 to 8,  y from 0 to 10
\put {$Figure \ 4.3.6$}  at 0 -1
\put {$\bullet$} at 0 0
\put {$\bullet$} at 0 8
\put {$\bullet$} at -2 8
\put {$e$} [t] at 0 -0.2
\put {$T'$} at  0  4
\put {$A$} at  0 8.7
\put {$B$} at -1 8.3
\put {$C$} at  -2 8.7
\put {$v_1$} [l] at 0.2 8
\put {$v_2$} [r] at -2.2 8
\put {$W_2$} at -2 10
\putrule from -3 9 to 1 9
\putrule from  -2 8 to 0 8
\setlinear \plot -4 4  0 0  4  4  0 8  -4 4  /
\setlinear \plot  -3 9  -2 8  -1 9  0 8  1 9  /
\setdashes
\setlinear \plot -6 6 -4 4  /
\setlinear \plot 4 4  6 6  /
\setlinear \plot -4 10  -3 9  /
\setlinear \plot -1 9  0 10  /
\endpicture
  \hss \egroup
\endinsert

Indeed, the arguments earlier applied to $W_1$ based at $v_1$ apply
just as well to $W_2$ based at $v_2$.

We know we can extend $W_2$ in a unique way to a sector
$\omega_2$ based at $e$.  We have no reason  to suppose
that $v_1$ will lie in $\omega_2$. However consider the vertex $v_3$
 of $T'$ in Figure 4.3.7.

\midinsert
  \line \bgroup \hss
\beginpicture
\setcoordinatesystem units <0.5cm, 0.866cm>
\setplotarea  x from -7 to 8,  y from 0 to 10
\put {$Figure \  4.3.7$} at 0 -1
\put {$\bullet$} at 0 0
\put {$\bullet$} at 0 8
\put {$\bullet$} at -2 8
\put {$\bullet$} at -1 7
\put {$e$} [t] at 0 -0.2
\put {$T'$} at  0  4
\put {$C$} at  -2 8.7
\put {$v_1$} [l] at 0.2 8
\put {$v_2$} [r] at -2.2 8
\put {$v_3$} [r] at -1.2 7
\put {$W_2$} at -2 10
\putrule from -3 9 to -1 9
\putrule from  -2 8 to 0 8
\putrule from -1 7  to  1 7
\setlinear \plot -4 4  0 0  4  4  0 8  -4 4  /
\setlinear \plot  -3 9  -2 8  -1 9  /
\setlinear \plot -2 8  -1 7  /
\setdashes
\setlinear \plot -6 6 -4 4  /
\setlinear \plot 4 4  6 6  /
\setlinear \plot -4 10  -3 9  /
\setlinear \plot -1 9  0 10  /
\endpicture
  \hss \egroup
\endinsert
\bigskip

 Since $v_3$ is the
unique vertex connected to $v_1$ and lying below it and to the left,
it must be the third vertex of the unique triangle containing and
lying below the edge from $v_1$ to $v_2$.  It is now clear that $v_3$
lies on a stretched gallery from $e$ to $v_2$, so we can be sure that
$v_3$ belongs to $\omega_2$.  It follows that $\omega_2$ contains
all of $T''$ ( Figure 4.3.8).

\midinsert
  \line \bgroup \hss
\beginpicture
\setcoordinatesystem units <0.5cm, 0.866cm>  
\setplotarea  x from -7 to 8,  y from -1 to 9
\put {$Figure \  4.3.8$} at 0 -1
\put {$\bullet$} at 0 0
\put {$\bullet$} at 0 8
\put {$\bullet$} at -0.5 7.5
\put {$e$} [t] at 0 -0.2
\put {$v_1$} [b] at 0 8.2
\put {$v_3$} [r] at -0.7 7.5
\put {$T$} at  0  2.5
\put {$T''$} at  -0.5 5
\putrule from -3.5 3.5 to 3.5 3.5
\putrule from 4 4 to 3 4
\putrule from -0.5 7.5 to 0.5 7.5
\setlinear \plot -4 4  0 0  4  4   /
\setlinear \plot -4 4  0 8  4 4    /
\setlinear \plot 3.5 3.5  -0.5 7.5  /
\endpicture
  \hss \egroup
\endinsert

Because we foresightedly added a layer
to our original $T$, it further follows that $\omega_2$ contains $T$,
i.e. $\omega_2\in\Omega(T)$.
One of the six boundary points of the apartment ${\cal A}$ is
represented by $W_2$, and so is equal to $\omega_2\in\Omega(T)$.
Similarly the boundary point represented by $W_3$ lies in $\Omega(T)$.
Since the other 4 boundary points were already known to lie in
$\Omega(T)$, we have found our desired apartment.\ep

\noindent {\bf Theorem 4.3.2} The action of $\G$ on $\O$ is topologically free. That is, for all
$g \in \G \backslash \{e\}$,
$$int\{ \o \in \O : g\o = \o \} = \emptyset .$$

\noindent {\bf Proof.} Suppose that  $g$ fixes all points of some nonempty open set
$V \subset \O$. By Proposition 4.3.1 , there is an apartment ${\cal A}$ all of whose boundary points
$\o_1, \dots , \o_6$ are contained in $V$.
 Now since the boundary points  $\o_1, \dots ,\o_6$  of ${\cal A}$ are fixed by $g$
it follows that ${\cal A}$ is stabilized by $g$. To see this,  fix an arbitrary vertex
$x_0$ in ${\cal A}$.  For $i=1,\dots,6$, the sectors  $S^{x_0}(\o_i)$
and $S^{gx_0}(g\o_i)$ are equivalent and hence contain a common subsector  $S_i$.
Thus  $S_i \subset {\cal A} \cap g{\cal A} \quad$ , $i=1,\dots,6$.  In particular
${\cal A} \cap g{\cal A} \ne \emptyset$. We may therefore suppose from the start that the
 base vertex $x_0$ lies in ${\cal A} \cap g{\cal A}$.  Now $g{\cal A}$ contains
$S^{gx_0}(g\o_i) = S^{gx_0}(\o_i)$, \quad $i=1,\dots,6$. But
 $\displaystyle{\cal A} = \cup_{i=1}^6 S^{gx_0}(\o_i)$. Therefore $g{\cal A} = {\cal A}$
and $g$ acts on ${\cal A}$ by some translation rather than some rotation or glide reflection.

Moreover all nearby apartments (i.e. all apartments having sufficiently many chambers in
common with ${\cal A}$) must be stabilized by $g$, since their boundary points will also be
in $V$ and hence fixed by $g$.  Translating by $g$ shows that if ${\cal A}$ and a nearby
apartment have one chamber in common, then they have infinitely many chambers in common.
But this is impossible, since there are apartments arbitrarily close to ${\cal A}$
having only finitely many chambers in common with ${\cal A}$.
To see this, simply choose a large convex subset $K$ of ${\cal A}$. Add a layer of triangles
to the boundary of $K$, none of which belongs to ${\cal A}$, and extend the resulting region
$K'$ to an apartment ${\cal A'}$. The contractibility of the building implies that
${\cal A} \cap {\cal A'} = K$.
\ep

\noindent 5. {\bf The algebraic structure of $C^*(s^{\pm})$}.

We can now show our main result.
\proclaim Theorem 5.1.
 The $C^*$-algebra  $C^*(s^{\pm}) = C(\Omega) \rtimes \Gamma$ is
simple and nuclear and canonically isomorphic to the reduced crossed product
$C(\Omega) \rtimes_r \Gamma$.

\noindent {\bf Proof}. Since the action of $\Gamma$ on $\Omega$ is minimal (Proposition
4.1.1) and topologically free (Proposition 4.3.2), it follows from [3: Corollary to Theorem 1] that
$C(\O) \rtimes_r \Gamma$ is simple. The remaining assertions follow from amenability of the
action  ( \S4.2).
\ep

\proclaim Corollary. Under the same hypotheses, the $C^*$-algebra $C^*_r(\Gamma)$ is subnuclear.

\noindent {\bf Proof}. $C^*_r(\Gamma)$ embeds in $C(\Omega) \rtimes_r \Gamma$.
\ep

\noindent {\bf Remarks}.1. It was shown in [7] that the group $\Gamma$ has Kazhdan's
property (T) and so is very far from being amenable.

2. A.M. Mantero and A. Zappa [18] have shown that, under the general assumptions of section 3.3,
the action of $\Gamma$ on $\Omega$ satisfies a certain geometric condition of [4]
which implies that $C^*_r(\Gamma)$ is simple.
The minimality of the action  (Proposition 4.1.1) and [4: Theorem 5] give another proof
that $C(\Omega) \rtimes_r \Gamma$ is simple .
\bigskip
\bigskip

\noindent 5.2 {\bf The linear case : $\widetilde A_2(F,v)$}

 Consider the class of affine buildings $\widetilde A_2(F,v)$
that arise in a natural way from linear groups [22: \S 9.2]. Not all the $\widetilde A_2$
buildings considered above can be constructed in this way. However for this class we
can obtain simplicity and nuclearity of the crossed product  using results of [21].

Let $F$ be a local field whose residual field has order $q$. There is a triangle
building $\triangle_F$ associated with $SL(3,F)$ [5:VI.9F],[22:\S 9.2].
Let ${\cal O}$ be the valuation ring of $F$ with respect to its valuation
$v$ and let $\pi \in F^{\times}$ be a uniformizer. i.e. $v(\pi)=1$.
A {\it lattice} in $F^3$ is a free ${\cal O}$-submodule of $F^3$ of rank 3.
Two such lattices $L$ and $L'$ are {\it equivalent} if $L'=aL$ for some $a \in K^{\times}$.
The vertices of ${\triangle_F}$ are the lattice classes. Its edges are pairs of vertices $x$
 and $y$ such that if $L$ is in the class of $x$ there is an $L'$ in the class of $y$
for which $\pi L \subset L' \subset L$. The chambers consist of triples of distinct
vertices $[L]$,$[L']$,$[L'']$ whose representatives can be chosen such that
 $L \supset L' \supset L'' \supset \pi L$.

The group $G=PGL(3,F)$ acts on ${\triangle_F}$ by defining $g\cdot[L]=[gL]$, for $g\in G$ and
each lattice $L$. A subgroup $\Gamma$ which acts simply transitively
 on the vertices of ${\triangle_F}$ in the case $F=F_q((X))$ is constructed explicitly in
[8:I \S4]. $G$ acts transitively on the boundary $\Omega$ and there is a point $\omega_0 \in
\Omega$ whose stabilizer is the group $P$ of upper triangular matrices in $G$ [5:Proposition
VI.9F]. Hence $\Omega$ is isomorphic to $G/P$ as a topological $G$-space.

 Assume that ${\triangle}={\triangle_F}$ and $G$,$P$,$\Gamma$
are as above.
   The group $P$ is solvable, hence amenable, and so the fact that
$C(\Omega) \rtimes \Gamma$ is nuclear and canonically isomorphic to $C(\O) \rtimes_r \Gamma$
follows from [21: Corollary 4.3].

 Moreover one can show that the action of $G$ on $\Omega$
 is topologically free directly as follows.

 Identify $\Omega$ with $G/P$ as above. Let $g \in \G \backslash \{e\}$ and
let $\omega=hP\in \Omega$
with $g\omega=\omega$  i.e. $h^{-1}gh=p\in P$.  We must find $\omega_0\in\Omega$ close to
$\omega$ such that $g\omega_0\ne\omega_0$. That is $h_0 \in G$ close to $h$ such that
$h_0^{-1}gh_0 \notin P$. Equivalently, we must find $k \in G$ close to $e$ such that
$k^{-1}pk \notin P$.

Now $p$ is upper triangular and $p\ne1$, so we may suppose that either $p_{12}\ne 0$
or $p_{11}\ne p_{22}$. Let $\varepsilon \in F$ and define $k\in G$ by
$$
(1-\varepsilon^2)k = \left [ \matrix
{1&\varepsilon&0 \cr
\varepsilon&1&0  \cr
0&0&1 \cr} \right]
$$
Direct calculation shows that the (2,1) entry of $k^{-1}pk$ is $\varepsilon(1-\varepsilon^2)^{-1}
(p_{22}-p_{11}-p_{12}\varepsilon)$, which is nonzero for suitable small $\varepsilon \in F$.
This proves the result.

Finally, the fact that  the action of $\Gamma$ on $G/P$ is minimal (\S 3.7) is analogous
 to the following
result [19: Section8]:  Let $G$ be a connected semisimple Lie group, $P$ a parabolic subgroup
of $G$ and $\Gamma_0$ a lattice in $G$. Then $\Gamma_0$ acts minimally on $G/P$.
\bigskip
\bigskip
\bigskip
\centerline{{\bf Appendix}:{\it Apartments and Galleries in Triangle Buildings}}

In this appendix we gather together various additional
definitions and results about triangle
buildings. The aim is to make it easier for someone  unfamiliar with buildings
 to read the later parts of the paper, particularly section 4.3. For further
information we refer to [5] and [22]

We deal exclusively with affine buildings of type $\widetilde A_2$. Such a building
is a simplicial complex $\triangle$ which is a union of certain subcomplexes called
${\it apartments}$. Each apartment $\Sigma$ is a Coxeter complex realized as a tiling of the
Euclidean plane by equilateral triangles. Any two simplices $A$,$B$ in the building
are contained in a common apartment. Furthermore, if $\Sigma$, $\Sigma'$ are apartments
containing $A$,$B$, then there exists an isomorphism between $\Sigma$ and $\Sigma'$
fixing $A$ and $B$ pointwise [5: Chapter IV]. We assume throughout that $\triangle$
has a ${\it complete}$ system of apartments: that is, the unique largest system of
apartments which contains every other system of apartments [5: Chapter IV.4]. The maximal
simplices in $\triangle$ are all equilateral triangles, and are called {\it chambers}.
 Two chambers $C$,$C'$ are {\it adjacent} if they have a common edge.
A {\it gallery} is a sequence of chambers
$$(C_0,\dots,C_d)$$
such that any two consecutive chambers $C_{i-1}$ and $C_i$ are adjacent. We
say that such a gallery connects $C_0$ to $C_d$.
 The integer
$d$ is called the length of the gallery. Any two chambers are certainly
 connected by a gallery since they lie in a common apartment.
 The distance between chambers
$C$ and $D$ in the building is the minimal length of a gallery $(C_0,\dots,C_d)$
connecting $C$ to $D$. A gallery which achieves this minimum is called a
{\it minimal gallery} from $C$ to $D$.
More generally, the distance between a chamber
$C$ and  a simplex $A$ in the building is the minimal length of a gallery $(C_0,\dots,C_d)$
with $C_0 = C$  and $C_d$ containing $A$. A gallery which achieves this minimum is called a
{\it stretched gallery} from $C$ to $A$.

Each vertex of $\triangle$ is labelled with a {\it type} which is an integer in $\{0,1,2\}$,
and each chamber has exactly one vertex of each type. The case we are interested in is
when the 1-skeleton of $\triangle$ is the Cayley graph of a group $\G$ with generating
set $P$. Recall from section 3.3 that the elements of $P$ correspond to the points in a
finite projective plane of order $q$. In this case
let $\tau : \Gamma \rightarrow \ZZ / 3 \ZZ$ be the homomorphism determined by $\tau
(a_x) = 1$ for each $x \in P$.  Then $\tau(g)$ is the type of $g$ for each $g \in \Gamma$.
Thus $e$ has type $0$, each generator $a_x$ has type $1$, and each inverse generator
$a_x^{-1}$ has type $2$.

Two adjacent chambers $C$,$C'$ are said to be {\it $i$-adjacent} if the vertex
of each which does not lie on their common edge has type $i$.
The type of a gallery $(C_0,\dots,C_d)$ is the sequence of labels $(s_1,\dots,s_d)$
such that  $C_{i-1}$ and $C_i$ are $i$-adjacent for $i = 1,\dots,d$. For example, in
Figure A.1 the gallery $(C_0,C_1,C_2,C_3,C_4,C_5)$ has type $(2,0,2,1,2)$.

\midinsert
  \line \bgroup \hss
\beginpicture
\setcoordinatesystem units  <1cm, 1.732cm>        
\setplotarea x from -2.5 to 5, y from -2 to 2    
\put {$Figure \  A.1$} at  0 -2
\put {$1$} [b] at -2  2.1
\put {$0$} [b] at  0  2.1
\put {$2$}    [b] at  2  2.1
\put {$0$} [l] at 3.1 1
\put {$1$} [l] at 4.1 0
\put {$0$} [t] at 3  -1.1
\put {$1$}  [t] at 1  -1.1
\put {$2$} [t] at  -1  -1.1
\put {$0$} [t] at -3 -1.1
\put {$2$} [r]  at  -4.1 0
\put {$0$} [r] at -3.2  1
\put {$2$} [l,b] at -0.7 1.1
\put {$1$}  [l,b] at 1.3  1.1
\put {$2$} [l] at  2.3 0.2
\put {$1$}  [b]  at  -2 0.2
\put {$0$} [t]  at  0  -0.2
\put {$C_3$} at -1 0.4
\put {$C_2$} at -1 -0.4
\put {$C_1$} at -2 -0.6
\put {$C_0$} at -3 -0.4
\put {$C_4$} at  0  0.6
\put {$C_5$} at  1  0.4
\putrule from -2   2     to  2  2
\putrule from  -3 1  to 3 1
\putrule from -4 0  to 4 0
\putrule from -3 -1  to  3 -1
\setlinear
\plot  -4 0  -2 2 /
\plot -3 -1  0 2 /
\plot -1 -1  2 2 /
\plot  1 -1  3 1 /
\plot  3 -1  4 0 /
\plot  -4 0  -3 -1 /
\plot  -3 1   -1 -1 /
\plot  -2 2   1 -1 /
\plot  0 2   3 -1 /
\plot  2 2   4 0 /
\endpicture
  \hss \egroup
\endinsert

It follows from [5: Chapter IV.3, Proposition 2] that a gallery is stretched
if it has the same type as a stretched gallery.

  Each apartment in the building $\triangle$ can be identified with a copy of the Euclidean plane
and hence inherits a natural metric which gives rise to a well-defined metric on the whole
building [5: Chapter IV.3]. A subcomplex  of $\triangle$
is called {\it convex} if it contains
every minimal gallery connecting any two chambers in it. The following result is very useful.
(See [5: Chapter VI.7 Theorems 1 and 2] for details.)

\noindent {\bf Theorem A.1}. A subcomplex $Y$ of an affine building of type $\widetilde A_2$
 which is either convex or
has nonempty interior and is isometric to a subset of the Euclidean plane $\RR^2$
is contained in an apartment.
Consequently $Y$ is an apartment if and only if $Y$ is isometric to $\RR^2$.

It is also important to note that $\triangle$ is topologically rather trivial.

\noindent {\bf Theorem A.2}. ([22: Appendix 4,page 185].)  An affine building is contractible.

\vfill \eject
\centerline{\bf References}.
\frenchspacing

\item{[1]} S. Adams, Boundary amenability for word hyperbolic groups and an application
to smooth dynamics of simple groups, {\it Topology} 33 (1994) 765-783.

\item{[2]} C. Anantharaman-Delaroche, Syst\`emes dynamiques non commutatifs et moyennabilit\'e,
Math. Ann. 279 (1987) 297-315.

\item{[3]} R. Archbold and J. Spielberg, Topologically free actions and ideals in
discrete dynamical systems,{\it Proc. Edinburgh Math. Soc.} 37 (1993) 119-124.
\medskip

\item{[4]} M. Bekka, M. Cowling and P. de la Harpe, Some groups whose reduced $C^*$-algebra
is simple, {\it Publ. Math. I.H.E.S.} (to appear).
\medskip

\item{[5]} K. Brown, {\it Buildings}, Springer-Verlag, New York 1989.
\medskip

\item{[6]} D. I. Cartwright, A. M. Mantero, T. Steger and
A. Zappa, Groups acting simply transitively on the vertices of a
building of type~$\widetilde A_2$, I and II,\ {\it Geom.\ Ded.}  {\bf 47}
(1993), 143--166 and 167--223.

\item{[7]}  D. I. Cartwright and W. M{\l}otkowski, Harmonic analysis for groups acting on triangle buildings, {\it J. Austral. Math. Soc. (A)} {\bf 56} (1994), 345--383.

\item{[8]}  D. I. Cartwright, W. M{\l}otkowski and T. Steger, Property (T) and
$\widetilde A_2$ groups, {\it Ann.\ Inst.\ Fourier} {\bf 44} (1993), 213--248.

\item{[9]} M.-D. Choi, A simple $C^*$-algebra generated by two finite order unitaries,
{\it Canad. J. Math.} 31 (1979) 867-880.
\medskip

\item{[10]} J. Cuntz and W. Krieger, A class of $C^*$-algebras and topological Markov chains,
{\it Inventiones Math.} 56 (1980) 251-268.
\medskip

\item{[11]} A. Fig\`a-Talamanca and C. Nebbia, {\it Harmonic Analysis and Representation Theory
 for Groups Acting on Homogeneous Trees.} Cambridge University Press,1991.
\medskip

\item{[12]} H. Furstenberg, Boundary theory and stochastic processes in homogeneous spaces,
in {\it Harmonic Analysis on Homogeneous Spaces.} (Symposia on Pure and Applied Math.)
Williamstown, Mass. 1972, Proceedings, vol.26 (1973) 193-229.
\medskip

\item{[13]} P. de la Harpe and G. Skandalis, Powers property and simple $C^*$-algebras,
{\it Math. Ann.} 273 (1986) 241-250.
\medskip

\item{[14]} P. de la Harpe, A.G. Robertson and A. Valette, Exactness of group $C^*$-algebras,
{\it Quart. J. Math.} (to appear).
\medskip

\item{[15]} G.Kuhn, Amenable actions and weak containment of certain representations
of discrete groups, {\it Proc. Amer. Math. Soc.} (to appear).

\item{[16]} G. Kuhn and T. Steger, More irreducible boundary representations of free groups,
{\it preprint}, University of Milan, 1992.

\item{[17]} A.M. Mantero and A. Zappa, Spherical functions and spectrum of the Laplacian
operators on buildings of rank 2,{\it Boll. U.M.I. (B)} 8 (1994), 419-475 .
\medskip

\item{[18]} A.M. Mantero and A. Zappa, The reduced group $C^*$-algebra of a triangle building,
 {\it preprint}, University of Genoa, 1993.

\item{[19]} G.D. Mostow, {\it Strong Rigidity of Locally Symmetric Spaces}, Ann. Math. Studies
No.78, Princeton University Press, 1973.
\medskip

\item{[20]} S. Mozes, Actions of Cartan subgroups,{\it preprint}, The Hebrew University,
 Jerusalem,1993.
\medskip

\item{[21]} J.C. Quigg and J. Spielberg,  Regularity and hyporegularity in $C^*$-dynamical
systems, {\it Houston J. Math.} 18 (1992), 139-152.
\medskip

\item{[22]} M. Ronan, {\it Lectures on Buildings}, Perspectives in Mathematics, vol.7,
Academic Press, 1989.
\medskip

\item{[23]} J. Spielberg, Free product groups, Cuntz-Krieger algebras, and covariant maps,
{\it International Journal of Maths.}, 2 (1991) 457-476.
\medskip
\bigskip
\bigskip
\bigskip
\bigskip

\noindent Department of Mathematics, University of Newcastle, N.S.W. 2308, Australia.\hfill\break
 [Electronic mail address: guyan@frey.newcastle.edu.au]
\bigskip
\noindent Department of Mathematics, Boyd Graduate Study Research Center,\hfill \break
 University of Georgia,
Athens, GA 30602-7403, U.S.A. \hfill\break
 [Electronic mail address: steger@margulis.math.uga.edu]

\end